\documentclass[a4paper,11pt]{article}
\usepackage[utf8]{inputenc}

\pdfoutput=1
\usepackage{amsmath}
\usepackage{amsfonts}
\usepackage{amssymb}
\usepackage{xcolor}
\usepackage[T1]{fontenc}
\usepackage{gensymb}
\usepackage{lmodern}
\usepackage{mdframed}
\usepackage{xcolor}
\usepackage[subtle]{savetrees}
\usepackage{enumitem,amssymb}
\usepackage{xurl}
\usepackage{hyperref}
\usepackage[hyphenbreaks]{breakurl}
\usepackage[center]{titlesec}
\usepackage[left=3cm,right=3cm,bottom=4cm]{geometry}
\usepackage{enumitem}
\usepackage{fancyhdr}
\usepackage{pdfpages}

\newlist{todolist}{itemize}{2}
\setlist[todolist]{label=$\square$}

\fancyheadoffset{0cm}

\fancypagestyle{plain}{%
  \fancyhf{}%
  \fancyhead[R]{\thepage}%
} 

\usepackage{titlesec}
\usepackage{sectsty}
\sectionfont{\color{black}}
\subsubsectionfont{\color{blue}}

\sectionfont{\centering}

\makeatletter
\renewcommand{\@seccntformat}[1]{%
  \ifcsname prefix@#1\endcsname
    \csname prefix@#1\endcsname
  \else
    \csname the#1\endcsname\quad
  \fi}
\newcommand\prefix@section{}
\makeatother
  
\def\valitemsep{-0.2em} 
\def\valrefleftindent{-2em} 
\def\valrefleftmargin{2em} 

\usepackage{etoolbox}
\patchcmd{\tableofcontents}{\contentsname}{\MakeUppercase\contentsname}{}{} 

\def\numberline#1{} 

\author{
Maurice Chiodo\footnote{Centre for the Study of Existential Risk, University of Cambridge, United Kingdom. \texttt{\href{mailto:mcc56@cam.ac.uk}{mcc56@cam.ac.uk}}.} \\
  Dennis M\"uller\footnote{RWTH Aachen University, Germany. \texttt{\href{mailto:dennis.mueller3@rwth-aachen.de}{dennis.mueller3@rwth-aachen.de}}.}    
}

\title{\vspace{0cm}Manifesto for the Responsible Development of Mathematical Works \\ \Large A Tool for Practitioners and for Management}

\begin{document}
\maketitle
\thispagestyle{empty}
\section*{SUMMARY}
This manifesto has been written as a practical tool and aid for anyone carrying out, managing or influencing mathematical work. It provides insight into how to undertake and develop mathematically-powered products and services in a safe and responsible way.  
\\Rather than give a framework of objectives to achieve, we instead introduce a process that can be integrated into the common ways in which mathematical products or services are created, from start to finish. This process helps address the various issues and problems that can arise for the product, the developers, the institution, and for wider society. 
\\To do this, we break down the typical procedure of mathematical development into 10 key stages; our ``10 pillars for responsible development'' which follow a somewhat chronological ordering of the steps, and associated challenges, that frequently occur in mathematical work. 
Together these 10 pillars cover issues of the entire lifecycle of a mathematical product or service, including the preparatory work required to responsibly start a project, central questions of good technical mathematics and data science, and issues of communication, deployment and follow-up maintenance specifically related to mathematical systems.
\\This manifesto, and the pillars within it, are the culmination of 7 years of work done by us as part of the Cambridge University Ethics in Mathematics Project. These are all tried-and-tested ideas, that we have presented and used in both academic and industrial environments. In our work, we have directly seen that mathematics can be an incredible tool for good in society, but also that without careful consideration it can cause immense harm. We hope that following this manifesto will empower its readers to reduce the risk of undesirable and unwanted consequences of their mathematical work. 

\newpage
\newgeometry{bottom=1cm}
\setcounter{tocdepth}{2} 
\tableofcontents
\restoregeometry
\newpage 
\section*{\uppercase{Introduction}}
\label{sec:intro}
\addcontentsline{toc}{section}{\nameref{sec:intro}}
\subsection*{Who this document is for:}
This manifesto was written as an aid to you, the practising mathematician, or indeed anyone doing or managing mathematically-powered work. It is very practical, and covers the fundamental and universal way in which a mathematical product or service is created, from start to finish, and the various issues and problems that can arise in the conception, production, use, and follow-up of mathematical work. We have not assumed that you have any training or background in ethics or responsible development, and we tried to keep our exposition as fundamental and elementary as possible, written in a way that we hope resonates with mathematicians as well as those in tangential areas such as computer science, physics, and engineering.  

Our intention here is to address the issues that we see as being common to all areas of mathematical work. As such, this tool is applicable and relevant to areas including machine learning, statistics, financial mathematics, cryptography, operations research, algorithms, and other forms of industrial and research mathematics. In the future, we hope to create tailored supplementary documents for each of the areas listed above, including some more domain and application specific documents. Supplementary material might contain worked through examples and explain in detail how the general manifesto applies to that specific area.

This is a tool for conscientious mathematicians who wish to avoid or minimise the harm that might come from their own direct mathematical work. But it is also a tool for highly-motivated mathematicians, or those in leadership positions, to assist other mathematicians around them in seeing the importance of, and moreover how to carry out, responsible mathematical work. We hope that the manifesto will help you in your own work, and aid you in helping others to improve their work as well.

\subsection*{How this document is laid out:}

The manifesto is laid out by first breaking down the typical procedure of mathematical development into 10 key stages; our ``10 pillars for responsible development''\footnote{A previous version of the pillars has been published as “Questions of Responsibility” in SIAM News: \url{https://sinews.siam.org/Details-Page/questions-of-responsibility-modelling-in-the-age-of-covid-19}.}. Each pillar addresses an important phase, and associated challenges, that frequently occurs in mathematical work. The ordering of the pillars has been carefully chosen, to reflect the order in which events, tasks, decisions, and actions generally happen in the mathematical development process. This is the \textit{approximate} order in which you might need to make decisions and take actions on the issues highlighted in each pillar. Think of this manifesto as a document that follows your workflow; as you move forward, the manifesto moves with you. But of course the mathematical development process is often nonlinear and iterative, and we therefore urge you to go back to previous pillars, or look ahead to later ones, when needed.

We have laid out each pillar with approximately 3 main ``steps'', to help you see how your work might evolve alongside the pillar. Within each step, we have included a handful of ``questions'' that you will likely need to address in order to manage the key idea of the step. These are deliberately phrased as questions, rather than instructions, because you will need to think about and reflect on them carefully. This manifesto is a tool, and not an automaton. It does not make the right decisions for you, but rather it helps you make the right enquiries and take the right decisions yourself. The questions are there to lead and guide you as you endeavour to create a responsible mathematical product.

To further assist you, we have broken down each question into more detailed ``checkboxes'', which are aspects that you should consider when trying to ensure that you have dealt with the question adequately. These checkboxes are not a check\textit{list} to be simply ticked off, and many of them will require substantial work to consider and address. Moreover, the checkboxes are not an exhaustive enumeration of things to address, but rather a minimal set, there to give you a starting point on some of the most common problems and pitfalls that stem from each question. At the end of each checkbox is a short layman’s reminder which we found useful when talking to other mathematicians; a few easy-to-remember words, written as a question, that encapsulate what the checkbox was trying to achieve.

You might notice that the early, and late, pillars, contain relatively little in the way of technical mathematical activity, compared to the middle pillars. This is to be expected, and as a mathematician you should still be acutely aware of these early and late pillars. There is a great deal of non-technical work to be done before you start doing heavy mathematics, and a great deal to be done after you complete most of your technical mathematics. Your expertise might be most closely aligned to the middle pillars, but your responsibility extends to both the early and late pillars.

We deliberately kept each pillar to 2 pages, to ensure this is a usable document. This is so you or a manager can easily bring a single pillar as one double-sided page to a team meeting or group discussion and say “Let’s look at this today”.

\subsection*{How to use this document:}
This tool is laid out in a linear fashion for clarity, but this can be deceiving. Going through this manifesto just at the start, or end, of your development work means you are treating ethical and responsible development as a bolt-on process rather than weaving it throughout your production and work as should be done.

To begin, you should go through and understand the whole manifesto before commencing your mathematical project, to identify how it initially fits in with your work plan and what initial actions you need to take. You then need to keep reviewing the entire document against your processes in an iterative way, going through all the pillars in detail. On each review pass, you will understand and see the relevance of a few more of the checkboxes. If you had total understanding of your work at all times (all the knowledge, all the resources, and the right framing of the problem) then you might be able to work through the pillars in a linear fashion. But since you're almost never in a perfect place with perfect insight, one pass of the manifesto is unlikely to suffice.

Different questions and checkboxes might have differing levels of relevance to your particular use of mathematics. It might be ok to skip some, but this should only be done with good reason. Not seeing the relevance of a question or checkbox is very different to establishing that it is not relevant, and as a mathematician you should appreciate this logical distinction. If you are unable to explain why a given question or checkbox is not relevant, then it might instead be wise to look more deeply to understand its possible relevance.

This manifesto will not only help you to \textit{solve} problems, but may also help you to \textit{detect} them. As you  try to address a problem that you identified from a particular checkbox, you might notice that the resolution of the problem may not be completely contained in the checkbox that made you aware of it in the first place. Deeper ethical issues usually arise because many things have gone wrong during the production process, and so it is natural that these problems may spill over to other checkboxes, too. You may need to go find, and act on, several checkboxes in order to properly address the problem. To put it into mathematical terms: there is no bijection between problems and checkboxes, it is a multi-valued map mapping one problem to potentially many checkboxes.

We have written this document with the assumption that your team consists of mathematically trained people with differing levels of ethical awareness\footnote{Chiodo, M. and Müller, D. (September 2018). Mathematicians and Ethical Engagement, \textit{SIAM News} 51, No. 9, p.6. Available online at \url{https://sinews.siam.org/Details-Page/mathematicians-and-ethical-engagement}.}, expertise and forms of training. Unlike many other ethical frameworks, this manifesto for responsible development cannot assume that its readership forms a homogenous group of professionals\footnote{Müller, D., Chiodo, M. and Franklin, J. (2022). A Hippocratic Oath for Mathematicians? Mapping the Landscape of Ethics in Mathematics, \textit{Science and Engineering Ethics} 28, article 41. Available online at \url{https://link.springer.com/article/10.1007/s11948-022-00389-y                                                          }.}. Thus, when using this document, you should be acutely aware of where your colleagues and collaborators come from. We have tried to include as many questions about your team, its awareness and the general human component of a project as possible, but we naturally could not cover it all. 

You will need to be aware that the people around you may understand a pillar, question or checkbox differently from you. This may concern the answer to a question, but it may also refer to a more fundamental understanding of the question or specific technical terms. Whenever possible, we attempt to use value-neutral and open language to cover as many mathematical disciplines as possible, but this will naturally lead to different interpretations among the readers. This puts the burden on you, the reader and user of this document, to talk through possible (mis)understandings with your peers when attempting to address individual questions and checkboxes. Even though it may seem tempting to rush some of the discussions or answers, in our experience, this can be damaging to the project's success and your team's integrity. It is critical to openly talk about some of the assumptions that everyone brings to the table at the very beginning when using this document and to potentially recapitulate them briefly before embarking on a new pillar.

\subsection*{Who produced this document:}
We are a team of mathematically-trained people, researchers and industry practitioners. Together, we have decades of experience studying and working in mathematics departments, teaching mathematics, conducting mathematical research, and working in various mathematically-powered industries. Alongside our other practical and mathematical experiences, over 7 years of our research in Ethics in Mathematics\footnote{Cambridge University Ethics in Mathematics Project. \url{https://www.ethics.maths.cam.ac.uk/}.} goes into this document. This manifesto is a culmination of all our experience doing mathematics, working with other mathematically trained experts and advising them in responsible development.

We previously published an early version of the pillars from this manifesto in SIAM News\footnote{Chiodo, M. and Müller, D. (September 2020). Questions of Responsibility: Modelling in the Age of Covid-19, \textit{SIAM News} 53, No. 7, p. 6-7. Available online at \url{https://sinews.siam.org/Details-Page/questions-of-responsibility-modelling-in-the-age-of-covid-19}.}, and we have presented and made use of versions of these pillars at a workshop on modelling covid organised by the Knowledge Transfer Network\footnote{\textit{Unlocking Higher Education Spaces - What Might Mathematics Tell Us?} Study group working paper (with M. Chiodo as a contributor), Turing Gateway to Mathematics, July 2020. Available online at \url{https://gateway.newton.ac.uk/event/tgmw82}.}, and other conferences, workshops and invited lectures, including at the London School of Economics\footnote{Chiodo, M. (February 2020). \textit{Guest lecture:  How do ethical issues arise in AI development and deployment?. London School of Economics}. MSc course on Algorithmic Techniques for Data Mining (MA429).}, the Royal Statistical Society\footnote{Chiodo, M. (2023). \textit{How to build an ethical blockchain (or any other mathematical technology).} Workshop on Responsible Modelling, hosted by the Royal Statistical Society. \url{https://rss.org.uk/training-events/events/events-2023/sections/responsible-modelling-(1)/}.} and QuantumBlack/McKinsey\footnote{Chiodo, M. (January 2020). \textit{Tech ethics - how hard can it be}. QuantumBlack / McKinsey \& Company, London.}. In addition, one of the authors has used various versions of these pillars in consultations with startups at Machine Intelligence Garage. These consultations were performed in conjunction with the existing MI Garage ethics framework\footnote{Machine Intelligence Garage (2021). \textit{Ethics Framework}. Available online at \url{https://migarage.digicatapult.org.uk/ethics/ethics-framework}.}; a document from which we have consequently also drawn some inspiration for this manifesto. When creating this general framework for responsible development in mathematics, we have further drawn inspiration from more technology-specific ethics frameworks, such as the pilot version of the \textit{Trustworthy AI Assessment List} ordered by the European Commission\footnote{AI HLEG (2018). \textit{Ethics Guidelines for Trustworthy AI.} Independent High-Level Expert Group on Artificial Intelligence, pp. 28-33. \url{https://digital-strategy.ec.europa.eu/en/library/ethics-guidelines-trustworthy-ai}.} and Microsoft's \textit{AI Fairness checklist}\footnote{Madaio, M.A. (2020). \textit{Co-Designing Checklists to Understand Organizational Challenges and Opportunities around Fairness in AI.} \url{https://www.microsoft.com/en-us/research/project/ai-fairness-checklist/}.}. 

Over the years, we have presented and worked through these pillars with students in our Cambridge seminar series on Ethics in Mathematics\footnote{Information about our regular Ethics in Mathematics seminar can be found on the websites of the Cambridge University Ethics in Mathematics Society (\url{https://cueims.soc.srcf.net/}) and of our research project (\url{https://www.ethics.maths.cam.ac.uk/)}.} and discussed them with many colleagues from various mathematical sciences, philosophy, and law, whose feedback we have hopefully successfully incorporated into this document. Many of the early and late pillars also draw ideas from discussions during one of the author's Fellowship at Auschwitz for the Study of Professional Ethics\footnote{Fellowships at Auschwitz for the Study of Professional Ethics: \url{https://www.faspe-ethics.org/design-technology/}}, the valuable and regular conversations on technology and ethics at the Centre for the Study of Existential Risk\footnote{Centre for the Study of Existential Risk, Cambridge. \url{https://www.cser.ac.uk/}} and lessons from \textit{EiM2: The second Meeting on Ethics in Mathematics}\footnote{Second meeting on Ethics in Mathematics (April 3-5, 2019). \url{https://ethics.maths.cam.ac.uk/EiM2/}.}. The ideas behind this manifesto are tried and tested in industry, research and teaching, and we hope that they will help you as much as they have helped us.

\subsection*{What this document is not:}
This manifesto is designed to help you identify the areas and issues that need to be addressed in your mathematical work, to assist you in developing your output responsibly and ethically. It does not give you direct answers; instead, it guides you to and gives you some of the questions to be answered. It is not an axiomatisation of responsible and ethical development, nor is it an algorithm to resolve these. Think of the work you do as being ``Turing-Complete'' metaphorically; the breadth of what you might produce is limitless, so such an ``ethics algorithm'' would be impossible for us to create.

This document is not a comprehensive enumeration of all possible questions that you need to ask of your product, or of yourself, in the development process. It covers several of the ones we see as important and relevant to most areas of mathematical work. But no finite document can cover every possible mathematical product ever to be made. This manifesto helps you see where to start, not how to finish.

And finally, this manifesto does not attempt to convince you of the importance of responsible and ethical development; we start with the assumption that this is something you already want to be doing. If you, or those around you, need further convincing of this, we recommend some of our other work on the importance of ethics in mathematics\footnote{Chiodo, M. and Clifton, T. (December 2019). The Importance of Ethics in Mathematics. \textit{LMS Newsletter} 484, 22-26. Available online at \url{https://www.lms.ac.uk/publications/lms-newsletter-back-issues}. Republished in the \textit{EMS Newsletter} 114, 34-37, December 2019.} and the value of using a code of ethics or a similarly codified document to maintain and foster responsible mathematical development\footnote{Müller, D., Chiodo, M. and Franklin, J. (2022). A Hippocratic Oath for Mathematicians? Mapping the Landscape of Ethics in Mathematics, \textit{Science and Engineering Ethics} 28, article 41. Available online at \url{https://link.springer.com/article/10.1007/s11948-022-00389-y    }.}. 

\newpage 
\section*{\uppercase{List of 10 key pillars for the responsible development of mathematical work}}
\label{sec:list}
{\Large
\begin{enumerate}
\item \textbf{Deciding whether to begin:} Why are you providing this mathematical product or service, and should you even do so?

\item \textbf{Diversity and perspectives:} Do your co-workers, superiors, and you have sufficient perspective, and do you understand the limitations and biases in your thinking?

\item \textbf{Handling data and information:} Are you using authorised and morally obtained datasets, in a responsible manner?

\item \textbf{Data manipulation and inference:} Do you have the expertise to properly manipulate data ensuring quality and ethics?

\item \textbf{The mathematisation of the problem:} What optimisation objectives and constraints have you chosen, and what are their real-life consequences? Who might the other impacted parties be?

\item \textbf{Communicating and documenting your work:} Are you properly considering how to comment and document your work and communicate the results to those who need them?

\item \textbf{Falsifiability and feedback loops:} Is your work falsifiable, and can you handle its large-scale impact and any feedback loops that arise?

\item \textbf{Explainable and safe mathematics:} Is your mathematical output explainable, and followed up with proper monitoring and maintenance?

\item \textbf{Mathematical artefacts have politics:} Are you aware of other non-mathematical aspects and the political nature of your work? What do you do to earn trust in yourself and your product?

\item \textbf{Emergency response strategies:} Do you have a non-technical response strategy for when things go wrong? Do you have a support network, including peers who support you and with whom you can talk freely?
\end{enumerate}
}
\newpage
\newgeometry{left=1.5cm,right=2cm,top=2cm,bottom=1cm}
\section{\uppercase{Pillar 1: Deciding whether to begin}}
\parbox{508pt}{\begin{center}{\Large \hspace{-4pt}Why are you providing this mathematical product or service, and should you even do so?} \end{center}}

\subsection[STEP 1: Solve the right problem]{\uppercase{Step 1: Solve the right problem}}
\subsubsection{Why do you want to solve this problem?}
\begin{todolist}
\setlength\itemsep{\valitemsep}
\item Consider why you are doing this project. Identify if your starting point is a problem worthy of being solved, or if you are simply seeking a use for your technology or mathematics. \textbf{Why are you doing this?}
\item Understand both the intended and possible use cases of the project’s output. Describe its intended goals and general purposes. Ensure that all goals are specific, measurable, achievable, relevant, and time-bound. Ensure that your goals and the (unintended) use cases of your product match. \textbf{What are you producing?}
\item Anticipate the costs to nature, and to society, of your project and its development process. Consider the impact if your project is successfully finished, and if it fails during or after development. \textbf{Can you anticipate the project's ethics?}
\item Clarify if the project as a whole has a beneficial impact. Pinpoint who is likely to benefit if your project is completed, and anticipate how those benefits are distributed. \textbf{What good will you do?}
\item Explore any possible general harms, including material cost, and any immaterial, hidden or indirect costs to society or nature. Pinpoint who will be negatively affected. \textbf{What harm might you do?}
\item Investigate the counterfactual, of what harm might occur if you don’t do your project, and also what good might come of you not doing it. \textbf{What happens if you don’t do this?}
\item Appreciate that the problems you want to solve deeply reflect who you are. Acknowledge that properly making the decision to build a product may require additional diversity in your team. \textbf{Do you have the breadth to even decide if this is worth doing?}
\end{todolist}

\subsubsection{How well do you understand your funders and beneficiaries?}
\begin{todolist}
\setlength\itemsep{\valitemsep}
\item Identify your funders, their motives, and what they can do with your work after you hand it over. Investigate their background and history for any questionable activities, or factors that may pose a risk to you or to others. Check who else they're giving money to, and investigate why they are funding you specifically. \textbf{Do you know who you’re dealing with?}
\item Check that your aims, objectives, and ideals are aligned with whoever is funding and influencing your work, to avoid future conflicts or compromising the project. \textbf{Are you and your funders aligned?}
\item Investigate who will be making use of your product, and what they might now be capable of doing with what you are providing them. \textbf{In whose hands are you voluntarily placing your solution or its parts?}
\item Recognise that overall control and direction of your product may pass to others on completion, possibly against your will, or in ways you disagree with. Establish who can end up with your product or its parts. Anticipate how much influence you will have over its eventual (mis)use, if any. \textbf{Who might unintentionally end up with final control of your product?}
\end{todolist}

\subsection[STEP 2: Solve the problem the right way]{\uppercase{Step 2: Solve the problem the right way}}
\subsubsection{
How are you weighing up the consequences of your solution?}

\begin{todolist}
\setlength\itemsep{\valitemsep}
\item Understand the distribution of benefits and harms more closely. Pinpoint where they will occur, and whether certain individuals, groups or institutions are more affected than others. Check if anyone is negatively affected at an existential level.  \textbf{Who wins, who loses?}
\item Establish agreed-upon quantitative and qualitative key performance indicators and other metrics that help you assess the project’s success, as well as the impact on those previously identified as positively or negatively affected. \textbf{How will you weigh up the pros and cons?}
\item Identify measures that introduce additional fairness and help to close gaps in society rather than perpetuating or exacerbating them. Adjust the project’s goals accordingly. \textbf{Can you improve fairness?} 
\end{todolist}

\subsubsection{Which issues with your mathematics are standard, obscure, or totally unexpected?}
\begin{todolist}
\setlength\itemsep{\valitemsep}
\item Review the common upsides and downsides of implementing and deploying mathematically-powered solutions to problems. \textbf{Have you checked all the usual mathematical issues?}
\item List further upsides, harms and implementation issues typical for your industry, your problem area, and for the mathematical tools you use. \textbf{How aware are you of your area’s pitfalls?}
\item Explore potential benefits and problems arising from your mathematics that are not well documented in the literature by consulting external domain experts. \textbf{Have you evaluated the obscure issues?}
\item Play devil’s advocate to understand how an irresponsible entity might utilise your mathematical output, tool or the project as a whole to deliberately inflict harm. Investigate if any sub-solutions and components of your product can be used for harm. \textbf{What deliberate harm is possible?}
\item Identify possible use case scenarios and inadvertent harms from your mathematics and how others might use it. Be specific about who will be harmed, and how. Check for any potential existential harm on an individual, group or societal level. \textbf{What unintentional harm is possible?}
\end{todolist}

\subsubsection{Why is a mathematical solution appropriate?}
\begin{todolist}
\setlength\itemsep{\valitemsep} 
\item Understand if using mathematics contributes positively to solving the problem. Identify whether introducing a mathematical solution might potentially cause harm. Avoid rationalising your use of mathematics just to be able to do more mathematics. \textbf{How does mathematics actually help here?}
\item Consider potentially less technical approaches and methods that don’t use the type, or complexity, of mathematics you are planning to implement. \textbf{Did you consider a non-mathematical approach?}
\item Establish whether a mathematical solution is the most appropriate approach. Investigate restricting the use of mathematics to certain subproblems. \textbf{Why is using mathematics the best option?}
\end{todolist}

\subsection[STEP 3: Allocate sufficient resources and responsibilities]{\uppercase{Step 3: Allocate sufficient resources and responsibilities}}
\subsubsection{Where does ethics sit within your project plan?}
\begin{todolist}
\setlength\itemsep{\valitemsep} 
\item Have a clear idea of how developing ethical mathematical products can be more costly in the beginning. Be prepared to up your game and appreciate the extra time, personnel, and expenditure constraints necessary to succeed in an ethical fashion. \textbf{Did you realise ethics takes effort?}
\item Allocate yourself and your team time to think about ethics. Understand what other project attributes you and your team will also have in mind during development, such as impact and profitability. Ensure people aren’t swamped with responsibilities and tasks. Realise that if everyone has to think about everything, then nobody is thinking anything fully through. \textbf{Where does ethics lie in your priority list?}
\item Know that without committing money and resources, your product is unlikely to be developed in an ethical way. Plan in the additional manpower, actions, time commitments, and training requirements needed to create an ethical product. \textbf{How will you budget for ethics?}
\item Avoid cutting corners or accumulating ethical debt, as an unethical project will be even more expensive in the long run. \textbf{Can you afford to carry out an unethical project?}
\item Prioritise and weigh up the ethics and project metrics. Consider whether achieving delivery and uptake of your project might overtake considerations surrounding its use and impact. \textbf{Will the interests of your project eclipse its purpose?}
\end{todolist}

\subsubsection{Is someone in charge, accountable and responsible at all times?}
\begin{todolist}
\setlength\itemsep{\valitemsep}
\item Follow good management practises by ensuring sufficient project oversight at a high level, and accountability of all parts at a low level. Identify and communicate those responsible for the project, its success, ethics, and all individual (non-)mathematical parts. \textbf{How has responsibility for all parts been assigned?}
\item Know and communicate to everyone who is managing and implementing responsible and ethical development throughout the team. \textbf{Who’s running your ethics programme?}
\item Propose mechanisms for testing and post-deployment monitoring that can quickly and effectively detect and respond to new problems. \textbf{Will there be any quality control?}
\item Contemplate what problems may arise out of the work, and what you might do if one does occur during the initial development or later life cycle. \textbf{Will you have a disaster plan?}
\end{todolist}
\newpage

\section{\uppercase{Pillar 2: Diversity and perspectives}}
\begin{center}{\Large Do your co-workers, superiors, and you have sufficient perspective, \\and do you understand the limitations and biases in your thinking?}    
\end{center}

\subsection[STEP 1: Setting up a diverse team]{\uppercase{Step 1: Setting up a diverse team}}

\subsubsection{Do you understand the value and importance of diversity?}
\begin{todolist}
\setlength\itemsep{\valitemsep}
    \item Appreciate that all good diversity work starts with introspection. Understand your own identity, biases, world view and perspective. \textbf{How well do you understand who you are?}
    \item Recognise the “problem-solving” value of having a diverse team of people working together, even if you cannot cover all perspectives. Understand that disagreements are necessary, and that you need perspective on the background to, solution methods for, and impact of, the problem. \textbf{Did you know that a diverse team produces a better output?}
    \item Honour the inherent social and moral value of diversity, beyond just the business case for it. Avoid using diverse people as a means to an end, but rather consider every person as an end in itself. \textbf{Do you value diversity beyond its profitability?}
    \item Be aware that problems such as bias in your product start long before you collect data; they start with problems and biases in your team and in your thinking. \textbf{How could your team’s biases manifest in your product?}
    \item Understand your own learning of diversity, and ask for help or further instruction if you are unsure that you are up to the task of incorporating diversity. Articulate what you’ve learned already. \textbf{What don’t you know about diversity?}
    \item Accept that diversity is not a one-time checkbox if you want to properly realise the benefits of having a diverse team. \textbf{Will you put in the effort to ensure diversity?}
    \item Understand that diversity is necessary at all levels of the project and institution to gain its full social, moral and economic benefits. \textbf{Will you diversify your tables of power?}
\end{todolist}

\subsubsection{What skills, viewpoints or diversity is your team currently missing?}
\begin{todolist}
\setlength\itemsep{\valitemsep}
    \item Identify the specific skills, training, ability, and expertise necessary for the success of your project, including non-mathematical skills such as relevant domain expertise or experience. \textbf{What knowledge beyond mathematics is needed in the room?}
    \item Establish diversity criteria for the project team, internal and external advisors. Aim for team members and advisors to have sufficient diversity of gender, ethnicity, age, education, etc. to properly solve the problem. \textbf{Do you have more than just diversity of skills and training?}
    \item Critically review the expertise, strengths and weaknesses of each team member over time. Ensure that potential weaknesses are resolved, and that additional training or hiring is done when needed, as the project may evolve and your needs may change. Have sufficient slack in your project's schedule and budget to allow for the maintenance of skills and diversity. \textbf{How are you maintaining an appropriate set of skills and expertise?}
    \item Critically review the background diversity of your team over time, beyond gender, age and nationality, as you need a diverse group of thinkers and problem-solvers. Be humble and appreciate what your team is missing. \textbf{How are you maintaining comprehensive representation in the team?}  
\end{todolist} 

\subsubsection{How will you create a team with a diverse background, expertise, and perspective?}
\begin{todolist}
\setlength\itemsep{\valitemsep}
    \item Identify what additional diversity you need to bring to the team and its management, before hiring. Write your job adverts and structure recruitment to foster it. Do all this for everyone you engage with to increase your understanding (incl. consultants, advisors and impacted parties). \textbf{Do you know who to bring in, and where to find them?}
    \item Consider, when hiring, what a diverse applicant is looking for. Appreciate the need for flexibility and understanding. \textbf{Do your job adverts make diverse people want to work with you?}
    \item Ensure the selection of new hires is done with diversity in mind. Weigh up the technical, non-technical, and social contributions a person brings to the team. \textbf{Are you selecting with diversity in mind?}
\end{todolist}
\subsection[STEP 2: Maintaining a diverse team]{\uppercase{Step 2: Maintaining a diverse team}}
\subsubsection{How will you bring a diverse team together to work?}
\begin{todolist}
\setlength\itemsep{\valitemsep}
    \item Understand that the project, your approach to the problem and your team’s culture need to be state-of-the-art to allow for diversity. Ensure everyone feels welcome, and can express their opinions, concerns and expertise freely in all interactions and discussions. Create sufficient safe spaces in your workspace, and communicate these to all team members. \textbf{Do your team’s values allow for diversity?}
    \item Champion trust, openness and accessibility. Put into writing how meetings are held to ensure fair treatment of all participants. Check that everyone understands their responsibilities and is accountable for their actions. \textbf{Does your workplace enable interactions that nurture diversity?}
    \item Maintain accessibility in all parts of the project, team, and work. Ensure that everyone can engage and participate fully in all interactions, information sharing, and decision making, and are not limited by physical, environmental, cultural, or technical barriers. \textbf{How is the workplace meeting diversity and accessibility requirements?}
\end{todolist}
\subsubsection{How will you harmonise a diverse team to get work done together?}
\begin{todolist}
\setlength\itemsep{\valitemsep}
    \item Reiterate to everyone the importance of talking with others about points and issues within the project. \textbf{Do team members know to interact with others?}
    \item Structure your work and teams so that diverse people are interacting. Ensure any partitioning of larger problems always respects the diversity requirements of each subproblem, including but not limited to expertise and experience. \textbf{Are people with different characteristics actually working together?}
    \item Establish ways for diverse people to talk and exchange ideas. Choose who chairs meetings, how problems are phrased, and what language is used. Appreciate that cross-disciplinary and cross-cultural collaboration is hard and is not solved by just gathering different people together. \textbf{Do diverse people have an effective way to interact?}
    \item Make sure everyone works together to see who may be harmed by the project. Understand that a diverse team is needed to identify all issues. \textbf{How will the team use their diversity to consider underprivileged groups?}
\end{todolist}

\subsubsection{Are you creating an environment which welcomes differing viewpoints?}
\begin{todolist}
\setlength\itemsep{\valitemsep}
    \item Prevent conflict, splintering, dominant majorities, or minority voices being dismissed. Ensure everyone is proactively listening to others. \textbf{How will you keep the peace and balance in a diverse team?}
    \item Establish ways to take on board the ideals and moral values of others in your team, in particular when these are expressed as such and not with the formalism of strict logical inference. Value everyone's voice, even when expressed differently to what you are used to. \textbf{How will you listen to people’s feelings?}
    \item Be aware that diverse teams have diverse needs. Constantly review and adjust the work, working conditions, and required accommodations. \textbf{What do you do when needs change?}
    \item Understand that human nature means you will need to invest effort in maintaining a diverse team and retaining members. Repeatedly check in with everyone and make sure that they are comfortable, seen and heard. \textbf{Do you actively check if everyone is ok?}
    \item Educate yourself and your team about groupthink. Ensure that it is safe to raise a differing opinion, and encourage it to the point of introducing a “devil’s advocate”. \textbf{Are you able to break groupthink?}
\end{todolist}
\subsubsection{How will you implement the output or recommendations of a diverse team?}
\begin{todolist}
\setlength\itemsep{\valitemsep}
    \item Make sure that all conversation is friendly and accepting, and takes on board the points of others. Have procedures in place for safe constructive feedback. Be prepared to up your game and defend the safe spaces of others. \textbf{Are you listening, enough to provide constructive feedback?}
    \item Grant those who may see problems the space and time to speak about them without fear of retribution. To improve your product and team, encourage discourse about ethical and social concerns of your product, its development and your team. \textbf{Do you reward people who find problems, or falsely punish them?}
    \item Be prepared to make changes after an issue has been raised. Maintain procedures that allow you to do so in a constructive and proactive way. \textbf{How are you making changes after hearing suggestions?}
\end{todolist}

\newpage
\section{\uppercase{Pillar 3: Handling data and information}}

\begin{center}{\Large Are you using authorised and morally obtained datasets, in a responsible manner? }    
\end{center}

\subsection[STEP 1: Training your team]{\uppercase{Step 1: Training your team}}
\subsubsection{How much do you respect data?}
\begin{todolist}
\setlength\itemsep{\valitemsep}
\item Ensure you have proper training on all standard techniques, tools, and methods for data collection and storage, including: privacy, copyright, database administration, informed consent, record-keeping requirements and the regulations pertaining to your problem domain. \textbf{Do you know your “data basics”?}
\item Understand that almost all data has an unavoidable human connection: either pertaining to the activities, behaviour, characteristics or preferences of people, or collected through systems built by people. \textbf{Can you understand the human connection of data?}
\end{todolist}
\subsection[STEP 2: Deciding what data to collect]{\uppercase{Step 2: Deciding what data to collect}}
\subsubsection{What data do you ideally want to collect and use?}
\begin{todolist}
\setlength\itemsep{\valitemsep}
\item Consider the data you need for the problem. Establish which parts of the data you already have, and which parts you still need. \textbf{What data do you need?}
\item Seek broad datasets that encompass all use cases of your product. Consider any weak points of your data, prior to collection. Search for any obvious red flags, such as missing demographics and edge cases for which you cannot collect data points. \textbf{Does your intended dataset cover enough to do the job?}
\item Avoid collecting an unwarranted quantity of data. Ensure that the volume of data is appropriate for the task at hand and the mathematics you want to use. \textbf{Are your datasets overkill for the task?}
\item Ascertain whether the data is actually helpful in solving the problem. Know that poorly selected, biased data may lead to harmful inferences. Pay particular attention to any data that may not be recent enough to positively contribute to your model or analysis. \textbf{Is all your data actually relevant?}
\end{todolist}

\subsubsection{How well do you understand the data you want?}
\begin{todolist}
\setlength\itemsep{\valitemsep}
\item Consider where you can feasibly source the data. Investigate the full data supply chain. Identify the true origins of your data: who generated it, where, how and why was it recorded, who it is about, and who the original legal and moral owner was. \textbf{Where and how does the data originate?}
\item Analyse the data to know what it contains, the insights it may reveal, and its deficiencies. Understand if you are working with structured or unstructured data, and the complications that come from it. \textbf{Do you know your data?}
\item Avoid treating data as “just data” and thus dehumanising it. Maintain awareness of the fact that, a lot of the time, your data is actually people’s data. \textbf{Does viewing data as an abstraction inhibit your empathy for its subjects?}
\item Identify any harm that could come to people if their data is not treated with respect. Prepare yourself well if you are to work with data about human characteristics, actions, activities, desires, and thoughts. \textbf{How might your use of data cause harm?}
\end{todolist}

\subsubsection{Is your intended data usage legal?}
\begin{todolist}
\setlength\itemsep{\valitemsep}
\item Understand the geography of the data you intend to use; where it is collected, processed, and used. Determine which legal jurisdictions are relevant, and the legal complications when crossing jurisdictions. Comply with all regulations, laws and best practices relevant to the data processing and your project. \textbf{Do you know the relevant data laws?}
\item Ensure that you have sufficient legal rights to utilise the data as you intend, including any limitations on profit generation, applications, and use cases. If necessary consult with your legal team or bring in external advice to check, and possibly to obtain, sufficient legal rights. \textbf{How are you allowed to use the data?}
\item Maintain supply chain due diligence. Whenever possible, check that all suppliers of your data collected and provided the data legally. Don't buy stolen, or otherwise illegally or immorally obtained data from someone else. \textbf{How does your data supply chain comply with the law?}
\end{todolist}
\subsection[STEP 3: Building a secure data environment]{\uppercase{Step 3: Building a secure data environment}}
\subsubsection{Do you understand the risks associated with your data?}
\begin{todolist}
\setlength\itemsep{\valitemsep}
\item Determine what risks are associated with your data, including privacy risks and data which is potentially interesting for malicious actors. Identify which data is of high, medium, or low risk, and protect it accordingly. Refrain from collecting any data that you cannot protect properly. \textbf{How risky is your data?}
\item Have in place a comprehensive emergency response plan for data loss (breach, leak, hack, etc). Consult appropriate experts in computer security to understand potential threats and data domain experts who can inform you about particularly vulnerable groups. \textbf{How have you prepared for data loss?}
\end{todolist}

\subsubsection{Do you have proper data privacy and transparency procedures in place?}
\begin{todolist}
\setlength\itemsep{\valitemsep}
\item Create a privacy policy, enforce it, and make it public where appropriate. Keep it up to date, reflecting both legal requirements and ongoing societal debates. Acknowledge that the legal requirements form the bare minimum of your privacy policy. \textbf{How is your privacy policy articulated and communicated?}
\item Understand when your data is no longer needed, and what forces and incentives may tempt you to keep data for too long. Prepare an auditable process for data to be securely erased, or further anonymised, when it or certain features are not required any more. \textbf{Do you practise good data hygiene?}
\item Design processes and easy-to-use interfaces that allow people to quickly detect and obtain any data you have on them, in accordance with any relevant laws. Ensure that your data discovery mechanisms do not pose any new privacy threats. \textbf{Can people find that you have their data?}
\item Design systems enabling people to excise their data from datasets, in accordance with any relevant laws. Whenever possible, create ways for people to excise their data from trained models, too. \textbf{Can people excise their data from you and your work?}
\item Map out what happens to all the data if the project concludes, or changes ownership. Execute appropriate data curation, and legal instruments, to preserve security and compliance of all data. Integrate this into your privacy policy, and properly communicate it. \textbf{Is your data and policies end-game robust?}

\end{todolist}

\subsubsection{Is your data properly secured?}
\begin{todolist}
\setlength\itemsep{\valitemsep}
\item Anonymise or pseudonymise your data, to preserve privacy and reduce the risk of losing sensitive data. Ensure that anonymisation happens during the collection or storage of the data, and not just at the level of database interfaces. Understand that interfaces are easy to remove and don't provide sufficient protection. \textbf{Are you “data minimising”?}
\item Properly encrypt data when possible. Ensure that access to the data is only given if absolutely essential. Brief all administrators about the data and its security. Ensure that administrators can question and challenge your team's data requests. \textbf{How are you controlling data storage and access?}
\item Keep your data protection measures up-to-date. Regularly inform yourself about recent threats and train your team accordingly. \textbf{Are you updating your data protection measures?}
\end{todolist}

\subsection[STEP 4: Collecting the data]{\uppercase{Step 4: Collecting the data}}
\subsubsection{Are you collecting and creating data in an ethical way?}
\begin{todolist}
\setlength\itemsep{\valitemsep}
\item Be aware that “legal” does not imply “moral”. Don’t obtain data through immoral means such as deception, subversions, or immoral scraping. \textbf{Are you tricking others into giving away their data?}
\item Ensure that people know the true value of their data. Avoid giving improper incentives to people in order to obtain their data, including any discounts that don’t match its true value, or exploitative monetisation of it. \textbf{Is your data acquisition ripping people off?}
\item Obtain full, proper, and voluntary informed consent for your intended data usage. Ensure this includes consent from the data creators and data owners, where possible. 
Be prepared to challenge any arguments telling you that it is economically or technically unfeasible to obtain consent for your data. \textbf{Do people consent to how you use their data?}
\item Avoid creating too much data through your project, product or service. Describe what will be created, and ensure this is communicated to the relevant users or participants. Make sure that you properly respect and protect all newly generated data. \textbf{Do people know what data you are collecting off them?}
\end{todolist}
\newpage
\section{\uppercase{Pillar 4: Data manipulation and inference}}
\begin{center}{\Large Do you have the expertise to properly manipulate data ensuring quality and ethics?}    
\end{center}

\subsection[STEP 1: Training your team]{\uppercase{Step 1: Training your team}}
\subsubsection{Are you providing sufficient training for your team to work with the data?}\begin{todolist}
\setlength\itemsep{\valitemsep}
\item Ensure that you train your technical team to understand, prevent and deal with common data problems (quality issues, data leakage, etc.). Educate the relevant non-technical staff about the common issues in an appropriate way in order to enable them to interact with the developers. \textbf{Is your team trained to manipulate data?}
\item Create a team awareness of the why’s and how’s of the data and information usage in the project. Promote discourse between the team covering which data is used, where it is used, and how it is used. \textbf{Is everyone aware of how data is used?}
\item Understand your impacted parties and what social biases or issues of discrimination they may experience. Acknowledge that this might be a good time to bring in additional domain expertise and a diverse group of people to look at your data. \textbf{Do you understand common biases for your problem domain?}
\item Train your team to detect biases in the data using both technical tools and hands-on understanding. Recognise that this may take a lot of resources if you or your team are doing it for the first time. \textbf{Does the team know how to detect biases in datasets?}
\end{todolist}

\subsection[STEP 2: Working with data while maintaining its quality]{\uppercase{Step 2: Working with data while maintaining its quality}}
\subsubsection{Are you working with relevant qualitative data?}
\begin{todolist}
\setlength\itemsep{\valitemsep}
\item Acknowledge that everyone can make simple mistakes when it comes to data. Investigate whether any of your data has been obtained from questionable sources, contains unintentional errors, or has been tampered with. \textbf{How are you checking for obvious red flags with the data?}
\item Acknowledge how the data was sampled or collected, and whether it sufficiently is current. Verify that your data distribution matches the application domain, and in particular that all relevant demographic groups are present in the data. Check if any groups are under- or overrepresented in the data, even in “hidden” ways such as having been the previous owners or sellers of the data. \textbf{Does the data distribution match your application domain?}
\item Understand if your data reflects a problematic world. Check for unwanted inequities in the systems your data relates to that your mathematical tool will likely perpetuate. Be prepared to intervene if you detect any. \textbf{Will your product perpetuate a problematic world?}
\end{todolist}

\subsubsection{Are you working with qualitative labelled data?}
\begin{todolist}
\setlength\itemsep{\valitemsep}
\item Select ethical features for the input X and and the output Y. Aim to maximise the coverage of your feature space whenever possible. Check if there are any dimensionality problems and ensure that your features work for your choice of algorithm. \textbf{Are you using good features?} 
\item Appreciate the human cost of labelling toxic data. Be prepared to provide sufficient psychological support for the data labellers. Ensure that they are appropriately reimbursed for their work, and any harms they suffer in the process of their work. \textbf{How are you caring for the labellers?}
\item Verify that your labelling is consistent, and that your human labellers can talk to each other. Understand any possible labelling ambiguities, and whether your labels capture uncertainty. Consider further technical and non-technical measures to improve the labelling process. \textbf{Are the labels of high quality?}
\item Understand that the quality of labels may decay in a static environment. Design sufficient improvement processes in your labelling system, such as having multiple humans labelling the same data points. \textbf{Will the labels get better over time?}
\item Ensure that the quality of your training data matches that of the production data. Check for circumstances that are different in the production environment, pay particular attention to those which could negatively impact the data quality. \textbf{Will the production environment surprise you?}
\item Understand that data pre-processing increases in complexity when you work with big unstructured data. Ensure that any data pre-processing does not negatively impact system performance. \textbf{Are you properly pre-processing the data?}
\item Understand that big data can have small data problems: understand outliers, periodic events, etc. \textbf{Are you addressing all the small data problems?}
\end{todolist}

\subsubsection{Are you properly documenting the problems of your data?}
\begin{todolist}
\setlength\itemsep{\valitemsep}
\item Chronicle the data and its processing from beginning to end. Record how it was gathered, what you did with it, and why you (pre-)processed it in a certain way. Highlight any non-standard collection and processing practices and everything that deviates from your well-established best practices. \textbf{How are you documenting the whole data pipeline?}
\item Understand what inferences you can and can’t make based on the data you have. Properly document all the limitations of the data, and any problems it may have, for future reference. Document the arguments and discussions that contributed to your understanding or lead to any relevant decisions. \textbf{How are you noting any unwanted properties of the data?} 
\item Understand and document what existing biases in the data your work could potentially eternalise, or even exacerbate. When possible use existing metrics for bias and fairness measurements, thus enabling your team to check if those biases have been counteracted. \textbf{Did you put all biases into writing?}
\end{todolist}
\subsection[STEP 3: Reuse the data whenever possible]{\uppercase{Step 3: Reuse the data whenever possible}}
\subsubsection{Is your data valuable to wider society?}
\begin{todolist}
\setlength\itemsep{\valitemsep}
\item Appreciate that data collection is resource intensive, and that others may not have access to such high-quality data. Understand the value of your data beyond its economic value. Acknowledge what good and bad the data you have curated brings to the world. \textbf{Are you aware of the value of the data?}
\item Think about what restrictions currently exist on the data, including ownership, privacy and fairness. Consider removing some of these restrictions if it is safe and permissible to do so, and leads to the data being more helpful for others and easier to share. \textbf{What further consent should you obtain to share the data more widely?}
\item Understand if your data can be used for good by others. See if sharing the data is legally and morally appropriate, and if so then consider sharing it. \textbf{Would it be helpful and safe for you to share the data?}
\item Decide with whom you would like to share the data, and what this may achieve. Investigate sharing with groups such as researchers, non-profits, government departments, or other businesses. Be prepared to discuss this with your team and any relevant stakeholders. \textbf{Do you want to share the data?}
\item Reflect on the value, restrictions and potential receivers of your data, and then deduce whether you are in a position to share it responsibly. Investigate if there is anything else preventing you from sharing the data or from making the relevant decision. \textbf{Can you effectively decide on whether to share the data?}
\end{todolist}
\subsubsection{How can you make the data available to others safely?}
\begin{todolist}
\setlength\itemsep{\valitemsep}
\item Perform a thorough background check of any entity to whom you supply data. Understand that due-diligence is not just about from whom you get the data, but also with whom you share data. Acknowledge that sharing data should follow similar legal and moral standards as collecting it. \textbf{Can you trust the entities you give data to?}
\item Think about implementing tools such as homomorphic encryption before allowing researchers or other teams to access your data for analysis. Restrict such work to an environment you control, rather than providing unencrypted, uncontrolled access. \textbf{Can entities make use of your data without taking a copy away?}
\item Obtain effective legal advice and instruments to safeguard the datasets you share. Implement additional mathematical tools, including data curation, anonymisation, and randomisation, for further assurance. \textbf{Are both the laws of the land, and of mathematics, protecting data you pass on?}
\end{todolist}
\newpage 
\section{\uppercase{Pillar 5: The mathematisation of the problem}}
\begin{center}{\Large What optimisation objectives and constraints have you chosen, and what are their \\real-life consequences? Who might the other impacted parties be?}    
\end{center}

\subsection[STEP 1: The bigger picture: Ethics of mathematisation]{\uppercase{Step 1: The bigger picture: Ethics of mathematisation}}

\subsubsection{Are you responsibly doing mathematics to solve the problem?}\begin{todolist}
\setlength\itemsep{\valitemsep}
\item Establish who has asked for or supported this work, and what they are trying to achieve. Make sure that you always have a responsible contact person for any questions that arise. \textbf{Who wants this, and why?}
\item Understand how your mathematics is used as part of a larger (corporate or political) strategy. Check whether your objectives and constraints match the bigger picture and potential future scenarios. Ensure you have a clear understanding of these goals, and are not just moving the goalposts to fit the best kick you can take. \textbf{How does your work fit into the ``bigger picture''?}
\item Identify all relevant entities (managers, product owners, external experts, colleagues, users and customers) who need to understand the technical choices and decisions. Be prepared to explain your decisions to your audience at an appropriate level. \textbf{Who needs to know and understand what you’re doing?}
\item Identify any situations in which people can defer responsibility. Ensure there is at least one person with a bird’s-eye overview of the project. Assign responsibility for checking that all assumptions are satisfied, that the mathematics and its associated processes work, and that any intermediate and final results are used properly. \textbf{Is someone taking full responsibility for the mathematics?}
\item Avoid conflicts of interest, including conflicts between your own personal interests and those of your impacted parties. Check for any hidden interests, and verify that you are in a good position to make decisions. \textbf{Does anyone have any conflicts of interest?}
\item Ensure that matters of plagiarism, copyright, and patent infringement, are properly monitored and avoided for the tools you use and build. Be aware of how such actions might harm others, as well as the legal and reputational ramifications for you and your product. \textbf{Are you stealing any ideas to produce ``your'' mathematical product?}
\end{todolist}

\subsubsection{Are you trying to do mathematics without ethics?}
\begin{todolist}
\setlength\itemsep{\valitemsep}
\item Ensure that every team member is aware of relevant ethical issues and actively considers these throughout their work and product life-cycle. \textbf{Is everyone thinking about ethics, regularly?}
\item Understand that you need to weave ethics into all of your work, and that there is no extended period where you can “just be thinking about mathematics”. Appreciate that this requires you to weave in more communication into the process of doing mathematics. \textbf{Is ethics continuously present in your thinking?}
\item Have a communication strategy so that everyone can raise questions and learn about ethics at any time. \textbf{Are you encouraging everyone, even technical people, to talk about ethics?}
\item Ensure that those engaging in mathematical work always have access to relevant other expertise (e.g., ethicists, project managers, legal advice), and vice versa. Avoid introducing any artificial segregation which silos teams or expertise. \textbf{Are tech and non-tech people communicating?}
\end{todolist}

\subsection[STEP 2: The ethics of your mathematics]{\uppercase{Step 2: The ethics of your mathematics}}
\subsubsection{How does your mathematical setup match the problem you want to solve? }\begin{todolist}
\setlength\itemsep{\valitemsep}
\item Analyse how you are converting physical and social phenomena into mathematical objects. Check that the world fits your mathematics, and isn’t being “forced” to fit through a myriad of abstractions and simplifications. Be aware of any additional assumptions or reductions that you introduce along the way. \textbf{How have you made the world fit into mathematics?}
\item Avoid contrived simplifications and don’t forsake correctness for mathematical beauty. Understand that the world is complicated and that overly simplified mathematics cannot properly address most edge cases. \textbf{Did you simplify the mathematics just to make  “nice”?}
\item Identify which aspects of the problem will be obscured through your abstraction and mathematisation. Be mindful of any social or physical aspects of the problem that cannot be measured or quantified effectively. \textbf{What do you lose in your mathematisation?}
\end{todolist}

\subsubsection{How are you using mathematical truth?}\begin{todolist}
\setlength\itemsep{\valitemsep}
\item Acknowledge that different uses and application domains might have different inherent levels of possible precision. Understand areas of highly imprecise, and precise, use of mathematics, and where your work lies on this spectrum. \textbf{How precise is the use case of your mathematics?}
\item Realise that academic mathematicians are often content with leaving elements out of their work, provided everything they say is ultimately true. Acknowledge that this approach may not serve you well, as mathematical truth is often meaningless without full context. Be aware that a practical solution is unlikely to be useable if you leave out what's trivial to you, or forget to communicate essential contextual knowledge, as what's easy for you can be hard for others. \textbf{Are you leaving the user to fill in the gaps?}
\item Understand the inspiration for your mathematics. Identify, for the problem you are working on, what (im)precise mathematics you are using, and how you are applying and interpreting your calculations. \textbf{Are you deploying precise mathematics meaningfully?}
\item Be aware that, while there are times when you can use mathematics to make very strong statements, these truths are often only relevant in a certain context or under certain assumptions, outside of which the (potentially) true statement may no longer be relevant. \textbf{Are you being misleading by saying something is provably true?}
\item Appreciate that your mathematics can still cause harm even when your mathematics gives provable, exact, reliable, and actionable outputs. Understand that true does not necessarily equate to good. Check if your true mathematics might be exploiting others and harming them. \textbf{Are you aware that true mathematics can still cause harm?}  
\end{todolist}

\subsubsection{How are you selecting what mathematics to use?}
\begin{todolist}
\setlength\itemsep{\valitemsep}
\item Understand your fundamental mathematical and development decisions and how deeply they determine your output. Justify why you are using a particular distribution, algorithm, model, simulation, tool, or approximation. \textbf{What mathematical shape did you assume the world takes?} 
\item Be conscious of all secondary mathematical decisions and their consequences, including boundary conditions, loss-functions, and assumptions. Check that you didn’t introduce elements that deviate from reality just to make the mathematics work. \textbf{What mathematical constraints did you impose on the world?}
\item Be aware of seemingly inconsequential mathematical decisions, in particular problems arising from algorithmic implementation. Be mindful of decisions that appear totally inconsequential, such a which programming language to use, but which actually heavily influence how your solution works. \textbf{What trivial mathematical decisions did you make without reflection?} 
\item Appreciate how initial systemic decisions affect your mathematics further down the line. Be prepared to question the decisions at each stage by developing a methodology for their justification and evaluation. Understand the consequences of these decisions separately, and when combined. \textbf{Do your mathematical decisions compound to deviate from reality?}
\end{todolist}

\subsubsection{Are you optimising for everyone's benefit? }\begin{todolist}
\setlength\itemsep{\valitemsep}
\item Understand what you interpret the optimal outcome of your work to be, and check if this is actually a sensible interpretation. Scrutinise any objectives provided to you by others. Understand how all this fits into the magic triangle of development constraints: quality, time and costs. \textbf{What are you optimising over?}
\item Identify the impacted parties beyond customers and users. Understand how they may be impacted through your work. Check who else might be using your mathematical product, and who they may be impacting with it. \textbf{Who does your mathematics affect?}
\item Work out how your mathematical solution projects into and affects the original physical or social problem. Understand how it is un-mathematised, interpreted, and used by you and others. \textbf{How is your mathematics projecting into the world?}
\item Account for direct and indirect benefits and harms to all relevant parties when deciding on objectives and carrying out your work. \textbf{Who will be helped, who will be harmed?}
\end{todolist}

\newpage
\section{\uppercase{Pillar 6: Communicating and documenting your work}}
\begin{center}{\Large Are you properly considering how to comment and document your work \\and communicate the results to those who need them?}    
\end{center}

\subsection[STEP 1: Internal communication]{\uppercase{Step 1: Internal communication}}
\subsubsection{Are you properly communicating and creating awareness about your mathematics?}\begin{todolist}
\setlength\itemsep{\valitemsep}
\item Appreciate and discuss the differing views held by your team and colleagues about ethics, mathematics and the overall project. Give everyone the room to talk about their viewpoint and reflect on others. Understand that you may need to have this conversation more than once to converge to a common understanding. \textbf{Is everyone on the same page?}
\item Arrange for everyone in your team or (indirectly) working with you to understand the assumptions underlying the project, its trajectory and mathematics. Pay particular attention to differences in understanding between the technical and non-technical parts of your team. \textbf{Does everyone know the project’s assumptions?}
\item Point out to all relevant colleagues the rationale behind your mathematical decisions, and inferences about data and information. Listen when your colleagues explain their mathematics to you. \textbf{How do you talk about your decisions with colleagues?}
\item Understand and internally communicate the transparency, fairness and reliability requirements for your mathematical work. Acknowledge the value of things not working out and talk about what does not work, too. \textbf{Do you talk about transparency and fairness?}
\end{todolist}

\subsubsection{Are you talking about the non-mathematical aspects of your work?}\begin{todolist}
\setlength\itemsep{\valitemsep}
\item Ensure that all responsibilities are properly communicated, and all areas of responsibility are established. Understand that this may be uninteresting for some of your colleagues, so do it in a way that keeps them engaged. \textbf{How are you talking about responsibilities?}
\item Make sure that the multitude of expertises and the different knowledges are reflected in a variety of communication means. Realise that not everything can be communicated via a \LaTeX \ document. \textbf{How are you matching the means of communication to fit others?}
\item Keep everyone updated on the project’s progress and problems at all times, including any changes in the project brief or application domain. Ensure that people can access past progress reports as this can help them to evaluate the overall trajectory of the project and their tasks. \textbf{Do you keep people updated?}
\end{todolist}

\subsubsection{How are you raising and documenting ethical issues about the mathematics?}\begin{todolist}
\setlength\itemsep{\valitemsep}
\item Ensure that there is a safe way to raise and document ethical issues specifically. Establish procedures that protect all ethics documentation from being modified improperly. Create guidance for safe meetings and for any minutes documenting the discussions. \textbf{Do you allow for ethical issues to be raised safely and securely?}
\item Have the documentation of ethical issues accessible to everyone, especially new hirees. Review and discuss it regularly. Be aware that your understanding of ethics develops as the project progresses. \textbf{Are you updating and talking about ethics regularly?}
\end{todolist}

\subsection[STEP 2: Internal documentation]{\uppercase{Step 2: Internal documentation}} 
\subsubsection{Are you documenting all your background decisions and reasoning?}\begin{todolist}
\setlength\itemsep{\valitemsep}
\item Comment all relevant mathematical or technical parts and decisions, such as governing equations, parameters, modelling, algorithms,  and source code, in a way that explains not only what was done, but why it was done. \textbf{Did you comment your decisions?}
\item Document any additional thoughts or discourse during development which may be necessary to understand or maintain the system in the future. \textbf{Did you comment your reasoning?}
\end{todolist}

\subsubsection{Are you continuously documenting your mathematical work?}\begin{todolist}
\setlength\itemsep{\valitemsep}
\item Establish procedures to properly and continuously document all of the mathematical work you produce. Set up criteria for when the documentation needs to be adjusted, such as when the mathematics, application, use cases or external circumstances change. Ensure that the documentation reflects the most recent version of your product. \textbf{How are you updating your documentation?}
\item Ensure that you meet the necessary legal record keeping requirements for your problem domain. \textbf{How are you recording the things that must legally be recorded?}
\item Make all documentation available to all relevant parties, such as any relevant mathematical experts, management, domain experts, users, and customers. Adjust the language and abstraction to fit the audience. \textbf{Can all relevant people access your documentation?}
\item Write readable documentation that allows the relevant parties, including those not part of the mathematical development, to analyse and audit the underlying mathematics, data usage and assumptions. Ensure that assumptions and limitations are made evident, and if necessary include additional explanations. \textbf{Can all relevant people understand your documentation?}
\item Prepare explanatory documents for any turnover in your team, as well as for your own benefit when reflecting back on work done earlier. \textbf{Are you making handover documents ensuring continuity?}
\end{todolist}

\subsection[STEP 3: Communication more widely]{\uppercase{Step 3: Communication more widely}}
\subsubsection{Are you speaking with domain experts?}\begin{todolist}
\setlength\itemsep{\valitemsep}
\item Understand the need to talk to external domain experts, including non-academically trained personnel. Acknowledge that users and customers can be experts, too. \textbf{Are you talking to domain experts?}
\item Adjust your communication and beware of different meanings, terminology, and assumptions. Learn where the domain experts come from and appreciate the unique guidance they can provide you with, even if their use of language differs from yours. \textbf{Are you and the domain experts speaking the same language?}
\item Ensure that both technical and non-technical people have access to domain experts, and allow for regular and repeated communication. Avoid just consulting them at the beginning or the end of the project. \textbf{Can everyone talk with the domain experts?}
\item Be intellectually humble and listen to the domain experts. Realise that domain experts don’t necessarily know the mathematics, but do have unique perspectives and insight knowledge that purely technical teams miss. \textbf{Are you listening to the domain experts?}
\end{todolist}

\subsubsection{How do you communicate with the external world?}\begin{todolist}
\setlength\itemsep{\valitemsep}
\item Endeavour to understand the mathematical knowledge, needs, fears and emotions of your user base, any other impacted parties, and the general public. Acknowledge that some of this may be hard to digest as it may contradict with your or your team's vision. Be prepared to talk through any challenges with your colleagues more than once. \textbf{How well do you know your audience?}
\item Establish exactly what you hope to achieve from such communications, what actions or changes you hope to see, and what you hope to learn. Whenever appropriate communicate it to those with whom you are engaging. \textbf{Why are you engaging in communication?}
\item Develop a communication method that is appreciated by others. Be aware that sometimes communication itself can have negative effects, even if it is true and accurate. Acknowledge that your technical training probably did not prepare you well for this task, and don't be shy to seek help with this whenever necessary. \textbf{How are you curating what you tell people?}
\item Communicate with the non-mathematical audience in a manner they comprehend. Realise that you need to tell them about more than just the mathematics or other technical decisions. \textbf{Are you effectively communicating with the public?}
\item Be aware of different communication speeds and the means of communication: communication needs outside of academia often require quick turnaround and succinct messaging. \textbf{Does your communication match the world around you?}
\item Ensure that all affected parties are made aware of any limitations, shortcoming, or risks of your mathematical product. Communicate these in a way that resonates with your audience. \textbf{Are you talking with people about what might go wrong?}
\end{todolist}
\newpage
\section{\uppercase{Pillar 7: Falsifiability and feedback loops}}
\begin{center}{\Large Is your work falsifiable, and can you handle its large-scale impact \\and any feedback loops that arise?}
\end{center}

\subsection[STEP 1: The falsifiability of your mathematics]{\uppercase{Step 1: The falsifiability of your mathematics}}
\subsubsection{What does failure mean for your mathematical system? }\begin{todolist}
\setlength\itemsep{\valitemsep}
\item Be aware of the kind of mathematics your system is built on. Understand that using data to create a model is different to hypothesising some properties of the system and setting up equations to solve based on that. \textbf{Is your mathematics data- or theory-driven?}
\item Learn how different mathematical systems typically fail, and how failure can be caused by reasons beyond faulty assumptions or wrong mathematics. \textbf{Do you understand why data- or theory-driven systems fail? }
\item Define, as precisely as possible, what an error or wrong output actually means for your product. Realise this goes beyond failing mathematical processes, and includes individual output scrutiny. \textbf{What does it mean for your particular product to be wrong?}
\item Identify what particular circumstances or ways exist for your mathematical product to bring up errors, be wrong, or fail. Include failed assumptions, edge cases, and wrong outputs. Note which of these error types are difficult to detect. \textbf{How can your mathematics fail?}
\item Understand if your mathematical predictions or computations can be falsified by a practical check, measurement, or comparison. Develop measures to test if your mathematics can be logically contradicted. \textbf{Is your work falsifiable in a measurable way?}
\item Determine if your mathematical approach has been scientifically verified via reproducible experimentation. Check if it can, or has, been shown to work robustly for your application. \textbf{Is the underlying mathematics scientifically proven for your application?}
\end{todolist}

\subsubsection{How is the application environment deviating from your initial assumptions over time?}\begin{todolist}
\setlength\itemsep{\valitemsep}
\item Identify if your application environment is changing, by looking at any periodic patterns, slow or rapid changes. Watch out for important aspects that might become invisible or unobservable over time. \textbf{Does your modelled reality change with time?}
\item Understand how individuals, groups or institutions are adjusting their behaviour, activities, and perceptions around your product. Consider all protected classes according to anti discrimination laws, as well as those who can’t access your product. \textbf{Are people changing to adapt to your product?}
\item Investigate what parts of the world your product might change, including markets, physical systems, and governance structures. Understand potential indirect changes, beyond its relationship with directly-impacted parties. \textbf{What else changes in the world when your product is deployed?}
\item Check if and how your mathematics is changing reality and what these changes are, particularly when it is used in social environments or as part of a larger socio-technical process. \textbf{How does your mathematics change the world?}
\item Quantify and measure how well your mathematics matches reality, flag when this deviation becomes significant, and note when such measurements are not possible. \textbf{Can you spot when your mathematics no longer matches reality?}
\item Identify changes in the application domain which may cause your system to fail. Be aware of potential failure from low probability events, shocks in the environment, updated regulations, competing systems, and slow unobservable environmental evolution. \textbf{What else can happen that will cause your mathematics to fail?} 
\end{todolist}

\subsection[STEP 2: The changing impact and its scale]{\uppercase{Step 2: The changing impact and its scale}}
\subsubsection{How much does your product impact the world?}\begin{todolist}
\setlength\itemsep{\valitemsep}
\item Understand how the impact pattern of your mathematics changes. Check for new users and impacted parties, beyond those identified originally. \textbf{Is the circle of impacted parties expanding?}
\item Identify an upper limit of impacted parties and entities that can be maintained in a sustainable fashion, paying particular attention to environmental and social concerns. Understand if, and how, your product might exceed your limit. \textbf{How big could this get?}
\item Check if a bug in your system, when deployed over a large application domain, might create a systematic weakness in society. Identify if your work changes the entire operation of a particular large ecosystem. \textbf{Is the impact scale compounding errors and creating systemic problems?}
\end{todolist}

\subsubsection{How does your work amplify ideas, patterns and problems?}\begin{todolist}
\setlength\itemsep{\valitemsep}
\item Understand that big data and big impact spaces are more valuable, useful, and harmful than the sum of their parts. Identify where your mathematics makes use of this. \textbf{Does your product leverage its scale and cause harm?}
\item Identify any existing social or physical problems relevant to your mathematical system, and how it might propagate or even amplify them. Appreciate that even a small development team can create a product that exacerbates social issues such as inequality. \textbf{Does your work potentially amplify existing problems?}
\item Weigh up whether your organisation’s operational capacity is sufficient to address the full impact of your work. Be aware that certain problems are too big to be undone. Appreciate the toll this may take on you professionally, financially, emotionally, and politically. \textbf{Are you prepared to deal with large-scale problems?}
\end{todolist}

\subsection[STEP 3: System interactions and feedback loops]{\uppercase{Step 3: System interactions and feedback loops}}
\subsubsection{How does your product interact with the world?}\begin{todolist}
\setlength\itemsep{\valitemsep}
\item Understand that your product might interact with other people, institutions and products that you didn’t expect, anticipate, or prepare for. \textbf{What else operates in your domain?}
\item Identify how any unexpected interactions of your product might create undesired outcomes. Monitor any changes in such interactions, and highlight which ones are only partially observable. \textbf{What interactions are there between your system and others?}
\item Trace how users and impacted parties are interacting with your product. Check if these interactions match what you had intended during the design and development phases. \textbf{How do users and your system interact?}
\item Understand if your mathematical system distorts peoples’ behaviour, emotions and actions, or the way they view and interpret various parts of the world around them. \textbf{Does your mathematics distort peoples’ perceived reality?}
\item Investigate and document ways people might misuse your mathematical product. Review and check it regularly against real usage. Understand that different forms of misuse can develop over time as people get to know your system better. \textbf{How might others (mis)use your work?}
\item Identify if your work incentivises unethical or harmful behaviour of individuals or institutions, even as side effects. Understand the larger social and market effects of your product. Establish ways to counteract such harmful behaviour. \textbf{What new behaviour does your product incentivise?}
\end{todolist}

\subsubsection{Can you identify and curtail specific feedback loops?}\begin{todolist}
\setlength\itemsep{\valitemsep}
\item Study the social or physical phenomena that your mathematics exploits and influences. Understand how the world changes as your mathematics is used. \textbf{How does your product change the world?}
\item Check if your product has any unwarranted dependencies on social or physical phenomena. Understand how your product inadvertently changes, or must change, to fit a changing world. \textbf{How does your product change as the world changes?}
\item Identify any higher order feedback loops that might occur as the world, and your product, keep changing together. Check if any are dangerous, and develop ways to stop them  escalating out of control. \textbf{How can you detect and dampen feedback loops?}
\end{todolist}
\newpage
\section{\uppercase{Pillar 8: Explainable and safe mathematics}}
\begin{center}{\Large Is your mathematical output explainable, and followed up with \\proper monitoring and maintenance?}
\end{center}

\subsection[STEP 1: Developing an explainable system]{\uppercase{Step 1: Developing an explainable system}}
\subsubsection{Do you know what you want your system to do?}\begin{todolist}
\setlength\itemsep{\valitemsep}
\item Articulate the outcomes and behaviour you expect or desire from your system, as well as that which you have specifically tried to avoid or prevent. Ensure that your system is properly aligned with the intended goals. \textbf{What do, and don’t, you want your system to do?}
\item Update your understanding of relevant legal requirements, codes of conduct, and standards, including any that are coming into existence. Include all rules and regulations that require you to prove your system works safely. \textbf{What are the rules, standards, and best practices of your industry?}
\item Evaluate the overall level of risk in your setting. Appreciate that this is the combination of intrinsic risk that comes from the problem domain itself, and extrinsic risk introduced through the solution you are creating. \textbf{How risky is your solution implementation?}
\end{todolist}

\subsubsection{How well do you understand what it does, and doesn’t, do?}\begin{todolist}
\setlength\itemsep{\valitemsep}
\item Ensure that your team has the expertise to scrutinise, change and repair your mathematical product at all times. This includes technical ability, but also domain or industry expertise. Allocate sufficient resources for any required additional training. Avoid having just one expert who understands a given part of your product. \textbf{Do you have the team expertise to scrutinise it?}
\item Perform a general safety check before deployment. Check that the mathematics, code and infrastructure are of sufficiently high quality. Test some outputs and operations to see they are sensible, useful, and safe. \textbf{How do you scrutinise your product pre-deployment?}
\item Verify which initial goals your system meets, and whether it is being used in unexpected ways. Check for instances of non-catastrophic failures which might nonetheless not meet ethical expectations. \textbf{Does the system’s use and behaviour meet your expectations?} 
\item Confine and restrict the behaviour of your system sufficiently to avert harm, and avoid outputs that are known to be uncertain. Place limiters on where your product can be deployed, and on how users can interact with it, to prevent unwanted outcomes. Investigate what failsafes you can put in to prevent harm in the event of system failure. \textbf{What limiters and failsafes do you put into the system?}
\item Understand how your system makes individual predictions or decisions, as well as its global behaviour. Identify what happens if the system is retrained or parameters are changed. Stress-test it to look for any bizarre, harmful behaviour. \textbf{Can it behave in an unexpected, dangerous way?}
\end{todolist}

\subsubsection{Can you explain its functionality and outputs?}\begin{todolist}
\setlength\itemsep{\valitemsep}
\item List to whom you need to explain your product’s workings. Understand what you need to explain, and what you hope to achieve through such explanations. \textbf{Who will you be explaining it to, and why?}
\item Establish which system properties you should explain, including overall operation and performance, average case behaviour, and individual outputs. Check for any laws and standards governing explainability and interpretability. \textbf{What functionality do you want to, or need to, explain?}
\item Identify which elements of your solution are more risky or dependent on unstable assumptions, and put in extra effort when it comes to interpreting and explaining these areas. \textbf{Is your mathematics sufficiently explainable, given the risks involved?}
\item Record where you used difficult to explain mathematics, particularly external solutions or uninterpretable statistical models. Identify the black and white boxes of your solution. Evaluate the trade-off between performance and interpretability of your product, choosing interpretable and explainable methods instead when the performance drop is not significant. \textbf{Are you using interpretable and explainable mathematics where possible?} 
\item Design in functionality enabling external audits of your system, especially if these are mandatory for your problem domain. Know what parts of the system can be audited, and what parts should be audited. Where possible, extend the audit to your team and processes. \textbf{How effectively can it be externally audited?}
\item Consider any relevant robustness constraints, where these come from, and how they might change in the future. Investigate if your system has any sensitive dependence on its formulation, training, or on computational inputs. Appreciate that robustness goes beyond robust mathematical techniques and includes a robust approach to interpretability and explainability. \textbf{Is your product robust and explainable in the long-run?}
\end{todolist}

\subsection[STEP 2: Monitoring, safety, and corrective steps]{\uppercase{Step 2: Monitoring, safety, and corrective steps}}
\subsubsection{How will you know when it is causing problems?}\begin{todolist}
\setlength\itemsep{\valitemsep}
\item Acknowledge that your system’s interaction with the world can never be fully understood during development. Appreciate the need to monitor its functionality and impact, potentially long after initial deployment. \textbf{How are you monitoring your system post-deployment?}
\item Instruct team members on where and how to report problems and issues they find, and create a supportive environment so that they can do so without fear of retribution. Implement a whistleblower policy to protect those who identify problems. \textbf{Can team members report problems, without retribution?}
\item Appreciate that problems, mistakes or misuse won’t always lead to obvious system failure. Realise that users and those affected by your product are sometimes best-placed to observe problems. Create structured, accessible ways for them to feed their observations back to you. \textbf{Can users or impacted parties report problems to you?}
\end{todolist}

\subsubsection{How are you keeping it safe after completion or deployment?}
\begin{todolist}
\setlength\itemsep{\valitemsep}
\item Understand the need for regular maintenance, quality assurance, and updates. Allocate enough materiel and personnel, and ensure they have the skills to detect problems and propose solutions. Grant them sufficient access and power to suggest or force changes. \textbf{Do you have an effective maintenance team?}
\item Distribute maintenance tasks to specific people or groups, and be clear on who is responsible and accountable for what. Ensure that all maintenance tasks are covered, but that nobody is overwhelmed. \textbf{Is responsibility, and accountability, effectively conveyed?} 
\item Determine how long after deployment you can, will, and need to, be able to perform changes and carry out maintenance. Follow the repairability standards of your problem area and legal jurisdiction. \textbf{How long after deployment will you maintain your system?}
\item Identify internal and external factors that may prevent you from doing updates. Aim to maintain sufficient system access and control after deployment to carry out changes and corrections effectively and quickly. \textbf{Can you make corrections, sufficiently quick?}
\end{todolist}

\subsubsection{What adjustments are you making to it based on your monitoring?}
\begin{todolist}
\setlength\itemsep{\valitemsep}
\item Use your monitoring output to analyse the problem at hand. Understand if it comes from a shift in the data, a faulty model or solution, or human error. Never rule out simple bugs before you have checked the system. \textbf{Can you tell if a human or technical corrective step is necessary?}
\item Accept that even in a non-emergency situation, having made a mistake or simply seeing your work perform badly can be challenging for anyone. Foster a supportive team culture, that is properly willing to make changes. \textbf{Is your work environment open for change?}
\item Realise that a local change to your system might solve one problem, but introduce several more. Check the functionality of the part you are changing, and any interdependencies with the rest of your system. \textbf{Are your adjustments creating more problems?}
\item Acknowledge that a big update can change or break the foundational decisions, peripheral choices, or ethics of the system. Understand that such updates may require you to revisit and rethink many of the key decisions you made earlier. \textbf{Does your update require a rethink of things?}
\item Understand the technical and mental challenges that come from having to take corrective steps. Appreciate that you or your colleagues may be heavily invested in those parts of the product that now need to change. \textbf{Are you caring for everyone whose work gets overwritten?}
\end{todolist}
\newpage
\section{\uppercase{Pillar 9: Mathematical artefacts have politics}}
\begin{center}{\Large Are you aware of other non-mathematical aspects and the political nature of your work? What do you do to earn trust in yourself and your product? }
\end{center}

\subsection[STEP 1: Identify the planetary and social costs]{\uppercase{Step 1: Identify the planetary and social costs}}
\subsubsection{Are you properly considering the environmental cost?}\begin{todolist}
\setlength\itemsep{\valitemsep}
\item Assess if your initial product has an acceptable ecological impact for what it does. Identify where else such impact occurs along your supply chain and development process. \textbf{Is your product and its development environmentally friendly?}
\item Ensure that you can deploy, maintain, and update in an environmentally friendly fashion. Understand the long-term environmental cost of your product. \textbf{Is your project environmentally sustainable?}
\item Understand how your product changes the consumption patterns of others. Attempt to predict what new enterprises your product might facilitate, and the potential negative environmental effects of these. \textbf{Is your project incentivising excessive consumption by others?}
\end{todolist}

\subsubsection{Are you properly considering the social cost?}
\begin{todolist}
\setlength\itemsep{\valitemsep}
\item Understand the ethics of labour requirements for mathematical products, such as the offshoring of programming or data labelling. \textbf{Are all workers being properly paid and treated well?}
\item Evaluate how much initial investment your project requires, and whether the economic return or social good validates this. Check where that money comes from, and what else could have been done with it instead. \textbf{Is your project a good use of money?}
\item Consider all foreseeable indirect costs, especially if your product replaces or competes with existing products, services or public infrastructure. \textbf{Who is hurt by what your product replaces?}
\end{todolist}

\subsection[STEP 2: Your mathematics has politics]{\uppercase{Step 2: Your mathematics has politics}}
\subsubsection{Do you understand the politics of your mathematics?}\begin{todolist}
\setlength\itemsep{\valitemsep}
\item Investigate whether there are entities, groups or individuals who are morally or politically opposed to your or your product, and why. \textbf{Who might your ideological enemies be?}
\item Reflect on how trustworthy and loyal your funders, supporters and friends are. Determine whether they will support you in tough times or when you are being challenged. \textbf{Who are your friends?}
\item Establish what your mathematical system changes or breaks. Understand those you might be directly competing with, as well as how any market, social, institutional or political structures will be negatively impacted by such changes. \textbf{Does your mathematics cause too much upheaval?} 
\item Check for any entities that may have both the desire and capability to harm or suppress your project or colleagues. Map out ways in which they could potentially do so. \textbf{Who might hurt you, and how?}
\end{todolist}

\subsubsection{Do you understand your politics?}\begin{todolist}
\setlength\itemsep{\valitemsep}
\item Understand your own political and ethical position in life, and how it fits into the larger environment you are working in. \textbf{What are your ethics and politics?}
\item Acknowledge that choosing when and how to solve a problem has inherent political meaning, even when your techniques are objective mathematics. Realise that you and your funder's values, political, social and moral understanding of the world define what problems you want to solve, and thus are never completely impartial. \textbf{Are you hiding from politics through mathematics?}
\item Check if there are individuals, groups or institutions that you or your product might interact with but who do not share your vision. \textbf{Whose toes are you treading on?}
\item Establish ways to evaluate who you would not work with on moral, political or otherwise ideological grounds. Ensure these are not purely emotive. Investigate and confirm any assumptions you make before deciding. \textbf{How will you decide who (not) to work with?}
\item Ensure your financiers, collaborators and business associates are vetted. Assess how you, your product, or impacted parties might be hurt by such associations. \textbf{How might your associations cause harm?}
\end{todolist}

\subsection[STEP 3: Earning trust from society]{\uppercase{Step 3: Earning trust from society}}
\subsubsection{How does society see your product?} 
\begin{todolist}
\setlength\itemsep{\valitemsep}
\item Describe in a candid way how your product was built and what information it draws upon, to clients and the public alike. Avoid hiding behind mathematical or technical language whenever possible. \textbf{How understandable is your whole operation?} 
\item Develop an honest marketing strategy, accessible and public product documentation and a structured way for the public to ask questions about your product and what data went into it. \textbf{Can everyone see what your product is, and does?}
\item Communicate widely the improvements and gains your work provides, as well as potential drawbacks and risks, to show that overall your product is a tool for good that benefits society. \textbf{Have you publicised the good that comes from your work?}
\item Identify how the public perceives your product, processes and team. Use this as an opportunity to improve your operations. \textbf{How well do you understand the public’s perceptions?}
\item Specify who could reject or be angered by your practices, product or vision. Understand how and why they might voice their opposition, and temper your product accordingly. \textbf{Will your solution be widely accepted by society?}
\end{todolist}

\subsubsection{Are you trustworthy?}\begin{todolist}
\setlength\itemsep{\valitemsep}
\item Acknowledge that trust is earned through clear and reasoned choices, responsible development, diligent maintenance, and demonstrable reliability. Avoid seeking blind faith in your or your product. \textbf{Do you know what it means to be trustworthy?}
\item Review all external solutions and tools, to maintain a trustworthy supply of third party components and mathematics. Concede that if you don’t understand the derivation of the mathematics you use, you cannot understand its limitations. \textbf{How trustworthy are your third-party components?}
\item Appreciate that the holistic combination of everything you do determines how trustworthy your system is. Be mindful of what mathematics you do, where and how you apply it, how you talk about it, the claims you make, your development processes, and how you and your team carry yourselves. \textbf{Are you and your mathematics trustworthy?}
\item Understand that using mathematics does not make your solution inherently trustworthy. Do not obfuscate deeper problems with mathematics. \textbf{Are you earning trust, or faking it through mathematics?} 
\item Avoid using mathematics and statistics to incentivise, force, or scare people into unwarranted actions or beliefs. \textbf{Are you misdirecting people through your mathematics?}
\end{todolist}

\subsubsection{Are you building a trustworthy product? }\begin{todolist}
\setlength\itemsep{\valitemsep}
\item Ensure that all subsystems of your product are trustworthy, not just output and interaction layers. Avoid using only bolt-on output monitoring or safety mechanisms, particularly because they can easily be removed later. \textbf{How is trustworthiness built into your system?}
\item Prepare manuals and other resources for users and customers that include the limits and dangers of the product. Articulate what your system can, and cannot, do. \textbf{Do people know your system’s capabilities and limitations?}
\item Appreciate that no product or system is perfect. Communicate, transparently, any possible failure points and the circumstances such as assumption breakdown that might lead to them. Realise that the key to trustworthiness includes knowing when not to use it or when it does not work reliably. \textbf{Can users understand when it fails, or know it has failed?}
\item Develop ways to manage expectations, so that people do not over-rely on your product. Ensure all your claims are justified. Prevent the deferral of responsibility to your product in inappropriate ways. \textbf{Are people placing undue faith in it?}
\item Understand that excess or unwarranted confidence in a product can be dangerous. Enable people to use it in its intended ways and aim to prevent any overuse or misuse. \textbf{Are people interacting with it safely?}
\item Explore how your product might, or could, be used in other domains beyond its original design. Articulate precisely where it is tested to work. Provide warnings about, and limiters against, using it in untested domains. \textbf{How might its use migrate to other domains?} 
\end{todolist}

\newpage
\section{\uppercase{Pillar 10: Emergency response strategies}}
\begin{center}{\Large Do you have a non-technical response strategy for when things go wrong? And a support network, including peers who support you and with whom you can talk freely?}
\end{center}

\subsection[STEP 1: Understand that things will go wrong]{\uppercase{Step 1: Understand that things will go wrong}}
\subsubsection{Do you understand the ways in which any project may go wrong?}\begin{todolist}
\setlength\itemsep{\valitemsep}
\item Understand project failure at different stages (during initial development, after deployment, in the long-run). Ensure sufficient safeguards are in place to protect impacted parties from harm, and yourself from liability claims. \textbf{What if your product collapses?}
\item Be aware that your product may exist for substantially longer than your team. Plan for when your team might disband, be it slowly by core members leaving, or abruptly when something goes wrong. \textbf{What if your team disbands?}
\item Understand that your ethical vision may differ from future stakeholders when your project changes ownership or institution. Ensure you’ve developed a robust team, processes and product where ethical considerations are built in, and not just bolt-ons. \textbf{What if you are acquired?}
\end{todolist}

\subsubsection{Do you understand the specific risks for your project?}\begin{todolist}
\setlength\itemsep{\valitemsep}
\item Understand the distribution of your risks, and potential impact. Prepare a list of emergencies with their core attributes that you can pass on to someone. Identify all the high probability, low risk-scenarios. Check for catastrophic low probability, high risk scenarios. \textbf{What can go wrong?}
\item Realise that the people you work with might have different objectives, and thus may not agree on what constitutes “wrong”. \textbf{Who is defining wrong?}
\item Analyse the best, most likely, and worst-case scenarios for each risk identified. \textbf{How bad can each problem get?}
\item Identify both the mathematical and non-mathematical response tools you need to resolve each problem. \textbf{What is needed to solve each problem?}
\item Understand the required response and resolution time for each problem. Check for interactions and resolution dependencies between various problems. \textbf{What is the timescale and interdependency of each problem?}
\end{todolist}

\subsection[STEP 2: Action plans for specific emergencies]{\uppercase{Step 2: Action plans for specific emergencies}}
\subsubsection{What are your emergency response capabilities? }\begin{todolist}
\setlength\itemsep{\valitemsep}
\item Develop mechanisms for internal and external people to raise the alarm. Admit if this is not possible for certain problems because of large time delays or insufficient observability of the problem space. \textbf{Is there a process for raising the alarm?}
\item Understand how quickly you can become aware of a problem, and how this compares to its urgency. \textbf{Can you spot problems fast enough?}
\item Figure out whether your processes and resources enable you to detect and react quickly enough to problems. Implement procedural changes if you find any disparities. \textbf{Can you respond to problems fast enough?}
\item Identify all mathematical and non-mathematical response tools at your disposal. Check if these are sufficient to deal with potential problems. \textbf{What is missing in your emergency response toolbox?}
\end{todolist}

\subsubsection{Do you have a non-mathematical response plan?}\begin{todolist}
\setlength\itemsep{\valitemsep}
\item Develop a clear plan of action for when a problem or emergency arises with your product. Incorporate non-mathematical actions into this, such as communications or physical activities and interventions. \textbf{Which parts of your response go beyond just doing more mathematics?}
\item Highlight who needs to do what in your response plan, and ensure that those with emergency roles have the resources, knowledge and motivation to deal with the fallout. \textbf{Have all emergency roles been assigned?}
\item Regularly update your understanding of potential disasters and best response strategies. Ensure your plan is properly communicated, understandable, and talked through on a regular basis. \textbf{Is your emergency plan clear, and regularly reviewed?}
\item Check that each team member understands their role in the emergency response plan, and rehearse it fully and regularly. Identify those not participating in the plan, and ensure they know what others are doing with a view to rendering assistance. \textbf{Does everyone know what to do in case of an emergency?}
\item Understand the need for individuals to have someone to talk to when things go wrong, for help and support. Realise that team members may need to talk through problems and their own individual circumstances, with co-workers and externally. Ensure that there are sufficient safe spaces, time and resources to do so. \textbf{Do you have someone to talk to when things get hard?}
\item Develop ethical strategies to protect yourself and your product, including legal responses, technological protection mechanisms and publicity from external attacks. \textbf{How can you defend yourself?}
\end{todolist}

\subsubsection{How do you execute a response plan?}\begin{todolist}
\setlength\itemsep{\valitemsep}
\item Identify the problem. Understand what happened, where and when it happened, and why it happened. Check who is affected, and how urgent the problem is. \textbf{What is the problem?}
\item Convene all relevant team members, and bring in any necessary external consultants or domain experts, as early as possible. Be aware that it is better to escalate early than not have sufficient resources. \textbf{How can you bring everyone into the room?}
\item Double check that the problem is sufficiently similar to your plan of action, and if not then quickly develop an appropriate plan from scratch. Evaluate what relevant resources you have at your disposal. \textbf{Does the problem fit your existing plan?}
\item Start executing your plan, and constantly monitor the situation. Adjust your plan if you identify any changes or unexpected hurdles. \textbf{Can you change your plan where necessary?}
\end{todolist}

\subsection[STEP 3: Debriefing and lessons learned]{\uppercase{Step 3: Debriefing and lessons learned}}
\subsubsection{Do you properly debrief after an emergency?}\begin{todolist}
\setlength\itemsep{\valitemsep}
\item Understand precisely what happened. Recap what part of your product or project failed, how you became aware of it, and how you dealt with it. \textbf{Do you know what happened?}
\item Pinpoint why the emergency occurred. Analyse how the combination of mathematical and non-mathematical factors contributed to it. Understand if the parts of your work that failed were within your control initially. \textbf{Why did the emergency happen?}
\item Identify what went right or wrong during your response. Avoid blaming people, remain objective, and be rational in your analysis. \textbf{What did or didn’t work in your response?}
\item Investigate why your response plan or parts of it did or didn’t work, and whether it was a successful response overall. \textbf{Why did your response succeed or fail?}
\end{todolist}

\subsubsection{Are you incorporating lessons learned and updating your emergency plan?}\begin{todolist}
\setlength\itemsep{\valitemsep}
\item Think about what can be learned from the emergency and its circumstances. Check if that particular emergency appeared in your initial plan, and if so, whether you correctly anticipated the risks and impact. \textbf{What did the emergency teach you?}
\item Think about what can be learned from the execution of your emergency plan; its successes and shortcomings. \textbf{What did your response teach you?}
\item Identify what to improve in your product, your personal and institutional politics, your team, and your emergency plan. \textbf{What will you change to be prepared for the future?}
\end{todolist}
\newpage
\restoregeometry
\section*{\MakeUppercase{Further Reading}}
\label{sec:reading}
\addcontentsline{toc}{section}{\nameref{sec:reading}}
\subsection*{{\color{blue}Pillar 1:}}
\begin{itemize}[label={},itemindent=\valrefleftindent,leftmargin=\valrefleftmargin]
\setlength\itemsep{\valitemsep}
\item EthicalOS (2018). \textit{Ethical Operating System Toolkit. A guide to anticipating the future impact of today's technology. Or: how not to regret the things you will build.} Available online at \url{https://ethicalos.org/}.
\item Schoemaker, P. (2015). The Power of Asking Pivotal Questions. \textit{MIT Sloan Management Review} 56(2), pp. 39-47. Available online at \url{https://sloanreview.mit.edu/article/the-power-of-asking-pivotal-questions/}. 
\item Spradlin, D. (2012). Are you Solving the Right Problem? \textit{Harvard Business Review.} Available online at \url{https://hbr.org/2012/09/are-you-solving-the-right-problem}.
\end{itemize}

\subsection*{{\color{blue}Pillar 2:}}
\begin{itemize}[label={},itemindent=\valrefleftindent,leftmargin=\valrefleftmargin]
\setlength\itemsep{\valitemsep}
\item Phillips, K. W., et al. (2014). How diversity works. \textit{Scientific American} 311(4), pp. 42-47. Available online at \url{https://doi.org/10.1038/scientificamerican1014-42        }.
\item Schiebinger, L. (2013). \textit{Gendered Innovations: How Gender Analysis Contributes to Research.} Report of the Expert Group “innovation through Gender”, pp. 1-144. Available online at \\ \url{https://genderedinnovations.stanford.edu/Gendered\%20Innovations.pdf}. 
\item Ely, R. J. \& Thomas, D. A. (2020). Getting Serious About Diversity. Enough Already with the Business Case. \textit{Harvard Business Review}, pp. 115-122. Available online at \url{https://hbr.org/2020/11/getting-serious-about-diversity-enough-already-with-the-business-case}. 
\end{itemize}

\subsection*{{\color{blue}Pillar 3:}}
\begin{itemize}[label={},itemindent=\valrefleftindent,leftmargin=\valrefleftmargin]
\setlength\itemsep{\valitemsep}
\item Davis, K. \& Patterson, D. (2012). \textit{Ethics of Big Data.} O’reilly Media, Inc, pp. 1-82.
\item European Commission (2021). \textit{Ethics and Data Protection}, pp. 1-22. Available online at \url{https://ec.europa.eu/info/funding-tenders/opportunities/docs/2021-2027/horizon/guidance/ethics-and-data-protection\_he\_en.pdf}.
\item Floridi, L. (2013). \textit{The Ethics of Information.} Oxford University Press, pp. 1-380. Available online at \url{https://doi.org/10.1093/acprof:oso/9780199641321.001.0001   }. 
\item Hand, D.J. (2018). Aspects of data ethics in a changing world: where are we now? \textit{Big Data} 6(3), pp. 176–190. Available online at \url{https://dx.doi.org/10.1089/big.2018.0083                }. 
\item Jo, E.S. \& Gebru, T. (2020). Lessons from archives: strategies for collecting sociocultural data in machine learning. \textit{FAT* '20: Proceedings of the 2020 Conference on Fairness, Accountability, and Transparency,} pp. 306-316. Available online at \url{https://arxiv.org/abs/1912.10389}.
\item Perez, C.C. (2019. \textit{Invisible Women.} Exposing Data Bias in a World Designed for Men. Chatto \& Windus, pp. 1-432.
\item Porter, M.A. (2022). \textit{A Non-Expert's Introduction to Data Ethics for Mathematicians.} arXiv Preprint:2201.07794. Available online at \url{https://doi.org/10.48550/arXiv.2201.07794    }.  
\item Wacks, R. (2015). \textit{Privacy: A Very Short Introduction.} Oxford University Press, pp. 1-176.
\end{itemize}

\subsection*{{\color{blue}Pillar 4:}}
\begin{itemize}[label={},itemindent=\valrefleftindent,leftmargin=\valrefleftmargin]
\setlength\itemsep{\valitemsep}
\item Gebru, T., et al. (2021). Datasheets for Datasets. \textit{Communications of the ACM} 64(12), pp. 86-92. Available online at \url{https://arxiv.org/abs/1803.09010}.
\item Gudivada, V. N., Apon, A., Ding, J. (2017). Data Quality Considerations for Big Data and Machine Learning: Going Beyond Data Cleaning and Transformations. \textit{International Journal on Advances in Software} 10(1/2), pp. 1-20. Available online at \url{https://www.researchgate.net/publication/318432363_Data_Quality_Considerations_for_Big_Data_and_Maachine_Learning_Going_Beyond_Data_Cleaning_and_Transformations}. 
\item Meyer, M. N. (2018). Practical Tips for Ethical Data Sharing. Advances in Methods and Practices in \textit{Psychological Science} 1(1), pp. 131-144. Available online at \url{https://dx.doi.org/10.1177/2515245917747656              }.  
\item Pipino, L.L., et al. (2006). \textit{Journey to Data Quality}. MIT Press. Available online at \url{https://doi.org/10.7551/mitpress/4037.001.0001   }. 
\item Reis, J. (2022). \textit{Fundamentals of Data Engineering: Plan and Build Robust Data Systems.} O’Reilly Media, pp. 1-422.
\end{itemize}

\subsection*{{\color{blue}Pillar 5:}}
\begin{itemize}[label={},itemindent=\valrefleftindent,leftmargin=\valrefleftmargin]
\setlength\itemsep{\valitemsep}
\item Eitzel M. V. (2021) A modeler's manifesto: Synthesizing modeling best practices with social science frameworks to support critical approaches to data science. \textit{Research Ideas and Outcomes} 7: e71553. Available online at \url{https://doi.org/10.3897/rio.7.e71553       }.
\item McKelvey, F. \& Neves, J. (2021). Introduction: optimization and its discontents. \textit{Review of Communication} 21(2), pp. 95-112. Available online at \url{https://doi.org/10.1080/15358593.2021.1936143              }.        
\item Passi, S. \& Barocas, S. (2019). Problem Formulation and Fairness. \textit{FAT* '19: Proceedings of the Conference on Fairness, Accountability, and Transparency}, pp. 39–48. Available online \url{https://doi.org/10.1145/3287560.3287567           }.
\item Selbst, A.D., et al. (2019). Fairness and Abstraction in Sociotechnical Systems. \textit{FAT* '19: Proceedings of the Conference on Fairness, Accountability, and Transparency}, pp. 59–68. Available online at \url{https://doi.org/10.1145/3287560.3287598            }.
\item Thompson, E. (2022). \textit{Escape from Model Land: How Mathematical Models Can Lead Us Astray and What We Can Do About It}. Basic Books, pp. 1-256. 
\item Williams, J. \& Gunn, L. (May 7, 2018). Math Can’t Solve Everything: Questions We Need To Be Asking Before Deciding an Algorithm is the Answer. \textit{Electronic Frontier Foundation.} Available online at \url{https://www.eff.org/deeplinks/2018/05/math-cant-solve-everything-questions-we-need-be-asking-deciding-algorithm-answer}.

\end{itemize}

\subsection*{{\color{blue}Pillar 6:}}
\begin{itemize}[label={},itemindent=\valrefleftindent,leftmargin=\valrefleftmargin]
\setlength\itemsep{\valitemsep}
\item Bhatti, J., et al. (2021). \textit{An Engineer’s Field Guide to Technical Writing}. Apress Berkeley, CA. 
\item Kearns, F. (2021). \textit{Getting to the Heart of Science Communication: A Guide to Effective Engagement}. IslandPress, Washington D.C., pp. 1-280. 
\item Mitchell, M., et al. (2019). Model Cards for Model Reporting. \textit{FAT* '19: Proceedings of the Conference on Fairness, Accountability, and Transparency}, pp. 220-229. Available online at \url{https://arxiv.org/abs/1810.03993}.
\item Parnas, D.L. (2011). \textit{Precise Documentation: The Key to Better Software}, pp. 125-148. In: Nanz, S. (eds) The Future of Software Engineering. Springer, Berlin, Heidelberg. Available online at \url{https://doi.org/10.1007/978-3-642-15187-3            \_8}. 
\end{itemize}

\subsection*{{\color{blue}Pillar 7:}}
\begin{itemize}[label={},itemindent=\valrefleftindent,leftmargin=\valrefleftmargin]
\setlength\itemsep{\valitemsep}
\item Ensign, D. et al. (2018). Runaway Feedback Loops in Predictive Policing. \textit{Proceedings of Machine Learning Research} 81, pp. 1–12. Available online at \url{https://arxiv.org/abs/1706.09847}.
\item Lum, K. \& Isaac, W. (2016). To predict and serve? \textit{Significance} 13 (5), pp. 14–19. Available online at \url{https://doi.org/10.1111/j.1740-9713.2016.00960.x                  }. 
\item Popper, C. (2022). \textit{The Logic of Scientific Discovery}. Routledge, pp. 1-544.
\item Stanford Encyclopedia of Philosophy (2022). The Problem of Induction. Available online at \url{https://plato.stanford.edu/entries/induction-problem/}. 
\item Jacobs, M. (2019). \textit{The Truth in the Falsification of Artificial Intelligence}, pp. 1-9. Available online at \url{https://soundideas.pugetsound.edu/cgi/viewcontent.cgi?article=1032&context=psupc}. 

\end{itemize}

\subsection*{{\color{blue}Pillar 8:}}
\begin{itemize}[label={},itemindent=\valrefleftindent,leftmargin=\valrefleftmargin]
\setlength\itemsep{\valitemsep}
\item Amodei, D., et al. (2016). \textit{Concrete Problems in AI Safety}. Available online at \url{https://arxiv.org/abs/1606.06565}.
\item Loyola-González, O. (2019). Black-Box vs. White-Box: Understanding Their Advantages and Weaknesses From a Practical Point of View. \textit{IEEE Access} 7, pp. 154096-154113. Available online at \url{https://doi.org/10.1109\%2FACCESS.2019.2949286}
\item Molnar, C. (2022). \textit{Interpretable Machine Learning. A Guide for Making Black Box Models Explainable (2nd Edition)}. Available online at \url{https://christophm.github.io/interpretable-ml-book}.
\item Phillips, P.J. (2021). Four Principles of Explainable Artificial Intelligence. \textit{National Institute of Standards and Technology Interagency or Internal Report} 8312, pp. 1-43. Available online at \url{https://doi.org/10.6028/NIST.IR.831      }.
\item Treveil, M. et al. (2020). \textit{Introducing MLOps. How to Scale Machine Learning in the Enterprise}. O’reilly Media, inc, pp. 1-183.

\end{itemize}

\subsection*{{\color{blue}Pillar 9:}}
\begin{itemize}[label={},itemindent=\valrefleftindent,leftmargin=\valrefleftmargin]
\setlength\itemsep{\valitemsep}
\item Botsman, R. (2017). \textit{Who Can You Trust?: How Technology Brought Us Together and Why It Might Drive Us Apart}. PublicAffairs, pp. 1-336.
\item Crawford, K. (2021). \textit{Atlas of AI: Power, Politics, and the Planetary Costs of Artificial Intelligence.} Yale University Press, pp. 1-336.
\item O'Neill, O. (2002). \textit{A Question of Trust}. Cambridge University Press, pp. 1-110.
\item O'Neill, O. (2018). Linking Trust to Trustworthiness. \textit{International Journal of Philosophical Studies} 26(2), pp. 293-300. Available online at \url{https://doi.org/10.1080/09672559.2018.1454637           }.
\item Rittberg, C.J. (2021). Intellectual humility in mathematics. \textit{Synthese} 199, pp. 5571–5601. Available online at \url{https://doi.org/10.1007/s11229-021-03037-3                }. 
\item Rogaway, P. (2015). \textit{The Moral Character of Cryptographic Work.} Available online at \url{https://web.cs.ucdavis.edu/~rogaway/papers/moral-fn.pdf}. 
\item Spiegelhalter, D. (2017). Trust in Numbers. \textit{Journal of the Royal Statistical Society, Series A} 180(4), pp. 949-965. Available online at \url{https://doi.org/10.1111/rssa.12302    }.
\item Winner, L. (1980). Do Artifacts Have Politics? \textit{Daedalus} 109(1), pp. 121-136. Available online at \url{https://www.jstor.org/stable/20024652}.

\end{itemize}
\subsection*{{\color{blue}Pillar 10:}}
\begin{itemize}[label={},itemindent=\valrefleftindent,leftmargin=\valrefleftmargin]
\setlength\itemsep{\valitemsep}
\item Matta, N. \& Ashhkenas, R. (2003). Why Good Projects Fail Anyway. \textit{Harvard Business Review}. Available online at \url{https://hbr.org/2003/09/why-good-projects-fail-anyway}.
\item Snedaker, S. and Rima, C.  (2014). \textit{Emergency Response and Recovery}, pp. 427-449. In: Business Continuity and Disaster Recovery Planning for IT Professionals (Second Edition). Available online at \url{https://doi.org/10.1016/B978-0-12-410526-3.00008-8           }.
\item Virilio, P. (2007). \textit{The Original Accident}. Polity Press, pp. 1-120.
\end{itemize}
\newpage
\addtocontents{toc}{\protect\enlargethispage{\baselineskip}} 
\thispagestyle{empty} 
\includepdf[pagecommand={\section[WALL POSTER OF THE 10 PILLARS ]{}}]{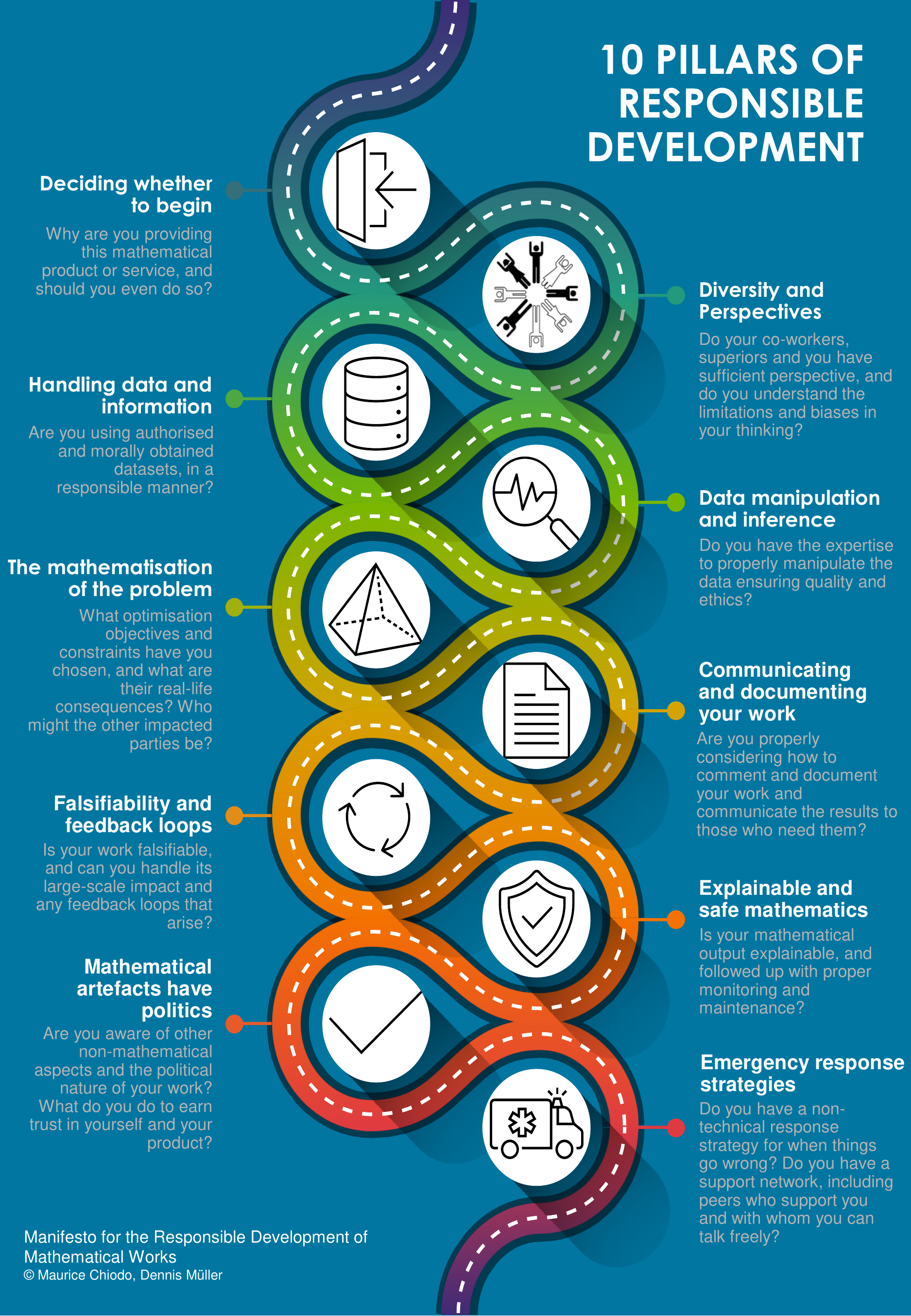} 
\end{document}